\documentclass[12pt]{amsart}
\usepackage[latin9]{inputenc}
\usepackage{bm}
\usepackage{amsbsy}
\usepackage{amstext}
\usepackage{amsthm}
\usepackage{amssymb}
\usepackage{xargs}[2008/03/08]
\usepackage[pdfusetitle,
 bookmarks=true,bookmarksnumbered=false,bookmarksopen=false,
 breaklinks=false,pdfborder={0 0 1},backref=page,colorlinks=false]
 {hyperref}
\hypersetup{
 pdfborderstyle=}

\makeatletter

\providecommand{\tabularnewline}{\\}

\numberwithin{equation}{section}
\numberwithin{figure}{section}
\theoremstyle{plain}
\newtheorem{thm}{\protect\theoremname}
\theoremstyle{remark}
\newtheorem{rem}[thm]{\protect\remarkname}
\theoremstyle{plain}
\newtheorem{lem}[thm]{\protect\lemmaname}
\theoremstyle{definition}
\newtheorem{defn}[thm]{\protect\definitionname}
\theoremstyle{plain}
\newtheorem{prop}[thm]{\protect\propositionname}
\theoremstyle{plain}
\newtheorem{cor}[thm]{\protect\corollaryname}
\theoremstyle{definition}
\newtheorem{example}[thm]{\protect\examplename}

\usepackage{amsfonts}
\usepackage{amscd}\usepackage{bm}
\usepackage{xcolor}\usepackage{xargs}[2008/03/08]

\numberwithin{equation}{section}

\definecolor{light-red}{rgb}{0.85,0,0}%
\usepackage{calrsfs}

\usepackage{tikz}

\numberwithin{figure}{section}
\numberwithin{table}{section}
\numberwithin{equation}{section}

\def\cP#1{{\cal P}_{#1}}
\def\cQ#1{{\cal Q}_{#1}}

\def\ss#1{\ifnum#1=1\let\next=\sone\else\let\next=\stwo\fi\next}
\def\tt#1{\ifnum#1=1\let\next=\tone\else\let\next=\ttwo\fi\next}
\def\ttm#1{\ifnum#1=2\let\next=\tone\else\let\next=\ttwo\fi\next}

\def\half#1{\frac{#1}2}
\def\hk{\frac k2}
\def\div{\mathord{\,|\,}}

\def\o{\omega}

\newdimen\rr
\newdimen\cc
\newdimen\dd

\def\olp(#1 #2 #3 #4){(\ifnum#1>1 #1\fi\ell,#2\ell\mid #3\ell,#4\ell)}
\def\tolp(#1 #2 #3 #4 #5 #6){(\ifnum#1>1 #1\fi\ell,#2\ell\mid #3\ell,#4\ell\mid\ifnum#5>1 #5\fi\ell,#6\ell)}

\newsavebox{\savepar}

\def\red{\color{light-red}}
\def\blue{\color{blue}}

\def\c#1{\ifcase#1\let\a=\relax\or
    \def\a{\frac k2-s} \or
    \def\a{\frac k2-v} \or
    \def\a{\frac k2+v} \or
    \def\a{\frac k2+s} \or
    \def\a{\frac k2-\half{s+v}-\half{t+u}} \or
    \def\a{\frac k2-\half{s+v}+\half{t+u}} \or
    \def\a{-\half{s+v}-\half{t-u}} \or
    \def\a{-\half{s+v}+\half{t-u}} \or
    \def\a{\frac k2+\half{s+v}-\half{t+u}} \or
    \def\a{\frac k2+\half{s+v}+\half{t+u}} \or
    \def\a{\half{s+v}-\half{t-u}} \or
    \def\a{\half{s+v}+\half{t-u}} \or
    \def\a{-\half{s+v}-\half{t-u}} \or
    \def\a{-\half{s+v}+\half{t-u}} \or
    \def\a{\frac k2+\half{s+v}-\half{t+u}} \or
    \def\a{\frac k2+\half{s+v}+\half{t+u}} \or
    \def\a{\half{s+v}-\half{t-u}} \or
    \def\a{\half{s+v}+\half{t-u}} \or
    \def\a{\frac k2+\half{s+v}-\half{t+u}} \or
    \def\a{\frac k2+\half{s+v}+\half{t+u}} \or
    \def\a{\frac k2+u} \or
    \def\a{\frac k2+t} \or
    \def\a{\frac k2-t} \or
    \def\a{\frac k2-u}
    \fi\vp\a}

\def\q#1{\ifcase#1\let\a=\relax\or
    \def\a{\frac k2-s}\or\def\a{\frac k2-v}\or\def\a{\frac k2+v}\or\def\a{\frac k2+s}\or
    \def\a{\frac k2-\half{s+v}-\half{t+u}}\or\def\a{\frac k2-\half{s+v}+\half{t+u}}\or
    \def\a{-\half{s+v}-\half{t-u}}\or\def\a{-\half{s+v}+\half{t-u}}\or
    \def\a{\frac k2+\half{s+v}-\half{t+u}}\or\def\a{\frac k2+\half{s+v}+\half{t+u}}\or
    \def\a{\half{s+v}-\half{t-u}}\or\def\a{\half{s+v}+\half{t-u}}\or
    \def\a{\half{s-v}-\half{t+u}}\or\def\a{\half{s-v}+\half{t+u}}\or
    \def\a{\frac k2+\half{s-v}-\half{t-u}}\or\def\a{\frac k2+\half{s-v}+\half{t-u}}\or
    \def\a{-\half{s-v}-\half{t+u}}\or\def\a{-\half{s-v}+\half{t+u}}\or
    \def\a{\frac k2-\half{s-v}-\half{t-u}}\or\def\a{\frac k2-\half{s-v}+\half{t-u}}\or
    \def\a{\frac k2+u}\or\def\a{\frac k2+t}\or\def\a{\frac k2-t}\or\def\a{\frac k2-u}
    \fi\a}
\def\cs#1{\ifcase#1\let\a=\relax\or
    \def\a{c_{1,3}}\or\def\a{c_{1,4}}\or\def\a{c_{2,3}}\or\def\a{c_{2,4}}\or
    \def\a{c_{1,5}}\or\def\a{c_{1,6}}\or\def\a{c_{1,7}}\or\def\a{c_{1,8}}\or
    \def\a{c_{2,5}}\or\def\a{c_{2,6}}\or\def\a{c_{2,7}}\or\def\a{c_{2,8}}\or
    \def\a{c_{3,5}}\or\def\a{c_{3,6}}\or\def\a{c_{3,7}}\or\def\a{c_{3,8}}\or
    \def\a{c_{4,5}}\or\def\a{c_{4,6}}\or\def\a{c_{4,7}}\or\def\a{c_{4,8}}\or
    \def\a{c_{5,7}}\or\def\a{c_{5,8}}\or\def\a{c_{6,7}}\or\def\a{c_{6,8}}\or
    \fi\a}
\def\ct#1{\ifcase#1\let\a=\relax\or
    \def\a{c_{4,2}}\or\def\a{c_{3,2}}\or\def\a{c_{4,1}}\or\def\a{c_{3,1}}\or
    \def\a{c_{6,2}}\or\def\a{c_{5,2}}\or\def\a{c_{8,2}}\or\def\a{c_{7,2}}\or
    \def\a{c_{6,1}}\or\def\a{c_{5,1}}\or\def\a{c_{8,1}}\or\def\a{c_{7,1}}\or
    \def\a{c_{6,4}}\or\def\a{c_{5,4}}\or\def\a{c_{8,4}}\or\def\a{c_{7,4}}\or
    \def\a{c_{6,3}}\or\def\a{c_{5,3}}\or\def\a{c_{8,3}}\or\def\a{c_{7,3}}\or
    \def\a{c_{8,6}}\or\def\a{c_{7,6}}\or\def\a{c_{8,5}}\or\def\a{c_{7,5}}\or
    \fi\a}

\newcommand{\p}[5]{\phi^{#1}\cdot \frac{\phi^{#2}-1}{\phi^{#3}-1}=\frac{\phi^{#4}-1}{\phi^{#5}-1}}
\def\pp#1(#2,#3,#4,#5){\p{#1}{#2}{#3}{#4}{#5}}

\let\cal=\mathcal

\def\bk{\mathbf{k}}

\renewcommand*{\backref}[1]{}
\renewcommand*{\backrefalt}[4]{{\tiny%
    \ifcase #1 {\red [Not cited.]}%
          \or [Cited on page~{\blue #2}.]%
          \else [Cited on pages {\blue #2}.]%
    \fi%
    }}
\renewcommand*{\backrefalt}[4]{}

\usetikzlibrary{calc}
\tikzset{
  octa/.pic={
     \foreach \i in {0,...,7}
       \draw (\i*45:0.65cm) node(\i){$\scriptstyle\bullet$};
     \node at (0,0) {$\scriptstyle\times$};
     \foreach \i in {0,...,7}
       \draw[dotted] let \n1={int(mod(\i+1,8))} in (\i.center) -- (\n1.center);
  },
  hexa/.pic={
     \foreach \i in {0,...,5}
       \draw (\i*60:0.65cm) node(\i){$\scriptstyle\bullet$};
     \node at (0,0) {$\scriptstyle\times$};
     \foreach \i in {0,...,5}
       \draw[dotted] let \n1={int(mod(\i+1,6))} in (\i.center) -- (\n1.center);
  }
}

\def\hhexa{
 \path (60:\x) pic {hexa}
	  -- (120:\x) pic {hexa}
	  -- (180:\x) pic {hexa}
	  -- (240:\x) pic {hexa}
	  -- (300:\x) pic {hexa}
	  -- (360:\x) pic {hexa};
}
\def\hhexagrapha{
     \foreach \i in {0,...,5}
       \draw (\i*60:\x) node(\i){$\scriptstyle\bullet$};
     \foreach \i in {0,...,5}
       \draw[dashed] let \n1={int(mod(\i+1,6))} in (\i.center) -- (\n1.center);
}

\def\hhexagraphb{
     \foreach \i in {0,...,5}
       \draw (\i*60:\x) node(\i){$\scriptstyle\bullet$};
     \foreach \i in {0,...,5}
       \draw let \n1={int(mod(\i+1,6))} in (\i.center) -- (\n1.center);
     \foreach \i in {0,...,5}
       \draw let \n1={int(mod(\i+2,6))} in (\i.center) -- (\n1.center);
     \foreach \i in {0,...,5}
       \draw[dashed] let \n1={int(mod(\i+3,6))} in (\i.center) -- (\n1.center);
}
\def\hexabasica{
     \foreach \i in {0,...,5}
       \draw (\i*60:\x) node(\i){$\scriptstyle\bullet$};
     \foreach \i in {0,...,5}
       \draw let \n1={int(mod(\i+1,6))} in (\i.center) -- (\n1.center);
}
\def\hexabasicb{
     \foreach \i in {0,...,5}
       \draw (\i*60:\x) node(\i){$\scriptstyle\bullet$};
     \foreach \i in {0,...,5}
       \draw let \n1={int(mod(\i+2,6))} in (\i.center) -- (\n1.center);
}
\def\hexabasicc{
     \foreach \i in {0,...,5}
       \draw (\i*60:\x) node(\i){$\scriptstyle\bullet$};
     \foreach \i in {0,...,5}
       \draw[dashed] let \n1={int(mod(\i+3,6))} in (\i.center) -- (\n1.center);
}

\makeatother

\providecommand{\corollaryname}{Corollary}
\providecommand{\definitionname}{Definition}
\providecommand{\examplename}{Example}
\providecommand{\lemmaname}{Lemma}
\providecommand{\propositionname}{Proposition}
\providecommand{\remarkname}{Remark}
\providecommand{\theoremname}{Theorem}

\begin{document}
\title{Overlaps in Field Generated Circular Planar Nearrings}
\author{Wen-Fong Ke and Hubert Kiechle}
\address{Department of Mathematics, National Cheng Kung University, Tainan
701, Taiwan}
\email{wfke@mail.ncku.edu.tw}
\address{Universität Hamburg, Fachbereich Mathematik, Bundesstr. 55, D-20146
Hamburg, Germany}
\email{hubert.kiechle@uni-hamburg.de}
\begin{abstract}
We investigate circular planar nearrings constructed from finite fields
as well the complex number field using a multiplicative subgroup of
order $k$, and characterize the overlaps of the basic graphs which
arise in the associated $2$-designs.
\end{abstract}

\maketitle

\section{Introduction}

\newcommandx\vp[1][usedefault, addprefix=\global, 1=]{\varphi^{#1}}%
\newcommandx\vpc[1][usedefault, addprefix=\global, 1=]{\phi^{#1}}%
\global\long\def\sone{j+s}%
\global\long\def\stwo{j-s}%
\global\long\def\tone{i+s}%
\global\long\def\ttwo{i-s}%
\global\long\def\bk{\boldsymbol{k}}%
\global\long\def\cT{\mathcal{T}}%
\global\long\def\cG{\mathcal{G}}%
\global\long\def\cI{\mathcal{I}}%
\global\long\def\cK{\mathcal{K}}%
\global\long\def\cN{\mathcal{N}}%
\global\long\def\cM{\mathcal{M}}%
\global\long\def\cO{\mathcal{O}}%
\global\long\def\cQk{\mathcal{Q}_{k}}%
\global\long\def\cC{\mathcal{C}}%
\global\long\def\cV{\mathcal{V}}%
\global\long\def\cE{\mathcal{E}}%
\global\long\def\bZ{\mathbb{Z}}%
\global\long\def\bC{\mathbb{C}}%
\global\long\def\bQ{\mathbb{Q}}%
\global\long\def\cPk{{\cal P}_{k}}%
\global\long\def\o{\omega}%
\global\long\def\Etot{\boldsymbol{\varphi}}%
\global\long\def\hk{\frac{k}{2}}%
\global\long\def\N{\mathbb{N}}%
\global\long\def\erc{E_{c}^{r}}%
\global\long\def\cP{{\cal P}}%
\global\long\def\cQ{{\cal Q}}%
\global\long\def\div{\mid}%
\global\long\def\vp{\varphi}%
\global\long\def\dD{{\cal D}}%
\global\long\def\Bp{{\cal B}_{\Phi}}%
\newcommandx\abs[1][usedefault, addprefix=\global, 1=]{\left|#1\right|}%
\newcommandx\gen[1][usedefault, addprefix=\global, 1=]{\langle#1\rangle}%
\newcommandx\norm[1][usedefault, addprefix=\global, 1=]{\left\Vert #1\right\Vert }%
\newcommandx\h[1][usedefault, addprefix=\global, 1=]{\frac{#1}{2}}%
\newcommandx\ttolp[6][usedefault, addprefix=\global, 1=]{(#1\ell,#2\ell\mid#3\ell,#4\ell\mid#5\ell,#6\ell)}%
\newcommandx\plist[5][usedefault, addprefix=\global, 1=]{(#4,#3\mid#1,#5)\text{ }\text{ with }\omega=#2}%
Planar nearrings were defined to connect nearrings and geometry. The
first three examples of planar nearrings were obtained by twisting
the multiplication of the complex number field, $\mathbb{C}$. These
examples provide many ideas for deriving geometrical and combinatorial
objects from planar nearrings. One example even inspires the notion
of circularity in planar nearrings (see~\cite{Clay92} for details).
In a circular planar nearring, ``circles'' are formed and one can
discuss the radii and centers of these circles just as one would with
circles in a complex plane. Understanding circular planar nearrings
also enables the creation of more circular planar nearring structures
in~$\mathbb{C}$.

One direction of research on circular planar nearrings involves selecting
equivalence classes, $\erc$, of circles with radius $r$ and centers
on another circle which has radius $c$ and center on $0$. Each $\erc$
has an associated graph, $G(\erc)$, which is naturally derived. This
graph is sometimes the union of spanning subgraphs called ``basic
graphs''. In other words, the ``overlapping'' of some basic graphs
produces the graph $G(\erc)$. In \cite{KeK96}, it is shown that
if a circular planar nearring $N$ is derived from a ring and $r$
is fixed, then the total number of basic graphs that appear in $G(\erc)$,
where $c\in N^{*}=N\setminus\{0\}$, is a function of $k$ only, even
if the nearring $N$ is changed to a different one. Since several
basic graphs can exist in a graph $G(\erc)$, the total number of
the graphs $G(\erc)$, where $c\in N^{\ensuremath{*}}$, varies from
one circular planar nearring to another.

In this work, we continue to study $\erc$'s for circular planar nearrings
constructed from finite fields as well the complex number field using
a multiplicative subgroup of order $k$. We begin with a brief review
of circular planar nearrings derived from fields, including some results
in \cite{KeK96}. We then define the overlaps of basic graphs and
show that for each $k$, there exists a finite set of primes, $\cQ_{k}$,
such that if $F=\text{GF}(q)$ is a Galois field of order $q$ with
$\text{char}F\not\in\cQ_{k}$ such that $k\mid(q-1)$, and $N$ is
a planar nearring constructed from $F$ using the multiplicative subgroup
of order~$k$ in~$F$, then the overlapping of the basic graphs
that occur in $N$ is exactly the same as that in $\mathbb{C}$ when
the regular polygon $C_{k}=\{z\in\mathbb{C}\mid z^{k}=1\}$ is used.

In section~\ref{S:normalized}, we discuss the normalized form for
overlaps which provide us the base to compare. In sections~\ref{sec:Main-Theorem}
and~\ref{sec:The-Proof-Main-Theorem}, with the help of a theorem
by Conway and Jones (Theorem~\ref{Th.CJ7}), we classify all overlaps
of basic graphs in $\mathbb{C}$. In the last two sections, we classify
all triple overlaps of basic graphs in $\mathbb{C}$, and conclude
that no further overlaps can be found.

The results we obtained in \cite{KeK96} have found applications (see
\cite{KeK95,KeK23-1,Kiechle94}). An application of the results obtained
in this paper to the number of solutions of equations $ax^{m}+by^{m}-cz^{m}=1$
over a finite field is in preparation.

\section{Preliminaries}

On a (left) nearring $(N,+,\cdot)$, the relation $=_{m}$ on $N$
given by $a=_{m}b$ if $ax=bx$ for all $x\in N$ is an equivalence
relation. When $N/{=_{m}}$ has at least three distinct classes and
if $ax=bx+c$ has a unique solution $x\in N$ for all $a,b,c\in N$
with $a\not=_{m}b$, we say that $N$ is \textit{planar}. Let $N$
be a planar nearring. For each $a\in N$, denote $l_{a}$ the map
from $N$ to $N$ given by $l_{a}(x)=ax$ for all $x\in N$. Then
the set $\Phi=\{l_{a}\mid a\in N$, $a\not=_{m}0\}$ is a fixed point
free automorphism group of the additive group $(N,+)$. The pair $(N,\Phi)$
is called the \textit{Ferrero pair }associated to~$N$.

Conversely, start with a Ferrero pair $(N,\Phi)$, where $(N,+)$
is a group and $\Phi$ a fixed point free automorphism group of $N$,
one can construct planar nearrings $(N,+,\cdot)$. For example, if
we take a field $F$ and a multiplicative subgroup $A$ of $F$ with
$|A|\geq3$, then $\Phi=\{l_{a}\mid a\in A\}$ is a fixed point free
automorphism group of $(F,+)$ and so $(F,\Phi)$ is a Ferrero pair.
We simply identify $\Phi$ with the set $A$ in this case. Any planar
nearring obtained from this Ferrero pair is referred to as a \textit{field
generated} planar nearring (see~\cite{Clay88} and~\cite{Clay92}).

Each planar nearring $N$ gives rise to certain combinatorial structures.
The one that concerns us here is an incidence structure. Let $(N,\Phi)$
be the Ferrero associated to $N$. For any $r,c\in N$, denote $\Phi r+c=\{\varphi(r)+c\mid\varphi\in\Phi\}$.
With ${\cal B}_{\Phi}=\{\Phi r+c\mid r,c\in N,r\not=0\}$, $(N,{\cal B}_{\Phi})$
is an incidence structure. If $N$ is finite, then $(N,{\cal B}_{\Phi})$
is actually a $2$-design (balanced incomplete block design).

In what follows, let $k\geq3$ be a fixed integer and $N$ a field
generated planar nearring with associated Ferrero $(F,\Phi)$ where
$(F,+,\cdot)$ is a field and $\Phi$ a multiplicative subgroup of
$F$ of order~$k$. Let~$\varphi$ be a generator of $\Phi$. We
will assume that
\begin{equation}
\big|(\Phi a+b)\cap\Phi c|\leq2\text{ for all }a,b,c\in F^{*}.\label{eq:circ}
\end{equation}
If this holds, $N$, as well as the Ferrero pair $(F,\Phi)$, is called
\textit{circular}.

Put $\bk=\{1,2,\dots,k-1\}$ and $\bk_{0}=\{0,1,2,\dots,k-1\}$, and
set
\[
\cI=\{(i,j,s,t)\in\bk^{4}\mid(i,j)\not=(s,t)\text{ and }(i,s)\not=(j,t)\}.
\]
Characterizations of circularity of $(F,\Phi)$ are given in \cite{Modisett89}
(see also \cite[\S5.3]{Clay92}).
\begin{thm}[{{{{{{\cite[Theorem~4]{Modisett89}}}}}}}]
\label{l:circ}The pair $(F,\Phi)$ is circular if and only $(\varphi^{i}-1)(\varphi^{t}-1)-(\varphi^{j}-1)(\varphi^{s}-1)\not=0$
for all \textup{$(i,j,s,t)\in\cI$}. This is equivalent to that $(\alpha-1)(\beta-1)-(\gamma-1)(\delta-1)\not=0$
for all $\alpha,\beta,\gamma,\delta\in\Phi\setminus\{1\}$ with $(\alpha,\beta)\not=(\gamma,\delta)$
and $(\alpha,\gamma)\not=(\beta,\delta)$.
\end{thm}

\begin{thm}[{{{{{{\cite[Theorem~8]{Modisett89}}}}}}}]
\label{thm:Pkfinit}For each integer $k>2$ there exists a finite
set of primes $\cPk$ such that for all finite fields $F$ and multiplicative
subgroup $\Phi$ of $F^{*}$ of order $k$, $(F,\Phi)$ is circular
if and only if $\textup{char}F\not\in\cPk$.
\end{thm}

The proof of Theorem~\ref{thm:Pkfinit} shows that $\cPk$ is the
union of the prime divisors of $k$ and those of the resultants $\text{Res}(g_{k},f_{i,j,s,t})$,
where $g_{k}=x^{k}-1$ and
\begin{equation}
f_{i,j,s,t}=(1-x^{i})(1-x^{t})-(1-x^{j})(1-x^{s}),\text{ }(i,j,s,t)\in\cI.\label{eq:fijst}
\end{equation}

Now $\Phi$ acts on $\Bp$ naturally: $\lambda\cdot(\Phi r+c)=\Phi r+\lambda c$
for all $\lambda\in\Phi$ and $\Phi r+c\in\Bp$. For $r,c\in F^{*}$,
denote by $\erc=\{\Phi r+\lambda c\mid\lambda\in\Phi\}$ the $\Phi$-orbit
of $\Phi r+c$ in $\Bp$. Then $|E_{c}^{r}|=|\Phi|=k$. It is known
that for any $r,r',c,c'\in F^{*}$, $\erc=E_{c'}^{r'}$ if and only
if $\Phi r'=\Phi r$ and $c'=\lambda c$ for some $\lambda\in\Phi$
(see \cite[(4.2)]{KeK96}). As $(F,\Phi)$ is circular, $\Phi r+c$,
$r,c\in F^{*}$ is regarded as a ``circle'' with radius $r$ centered
at $c$, and hence $\erc$ is the family of circles of radius $r$
centered at the points of the circle $\Phi c=\{\lambda c\mid\lambda\in\Phi\}$.

To visualize $\erc$ consider $(\mathbb{C},U)$ where $U=\{z\in\mathbb{C}\mid z^{k}=1\}$,
the regular $k$-gon inscribed in the unit circle $C=\{z\in\mathbb{C}\mid|z|=1\}$.
Thus, ${\cal B}_{U}$ is the collection of all regular $k$-gons in
the complex plane, and for $r,c\in\mathbb{C}^{*}$, $\erc$ is the
collection of the $k$ regular $k$-gons with radius $|r|$, centered
at $\lambda c$, $\lambda\in U$.
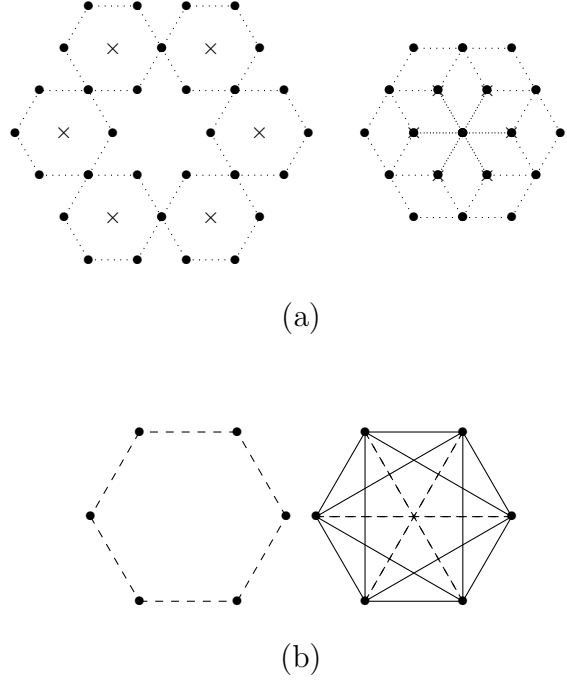
\begin{figure}
\begin{centering}
\begin{minipage}[b][1\totalheight][t]{10cm}%
\begin{center}
\begin{tikzpicture}
\begin{scope}
\newcommand\x{
1.3cm}
\hhexa
\end{scope}
\begin{scope}[xshift=4cm]
\newcommand\x{
0.65cm}
\hhexa
\end{scope}
\end{tikzpicture}\hspace*{0.35cm}
\par\end{center}

\begin{center}
(a)
\par\end{center}
\begin{center}
\vspace*{0.5cm}
\par\end{center}%
\end{minipage}\\
\begin{minipage}[b][1\totalheight][t]{8cm}%
\begin{center}
\begin{tikzpicture}
\newcommand\x{
1.3cm}
\begin{scope}
\hhexagrapha
\end{scope}\begin{scope}[xshift=3cm]
\hhexagraphb
\end{scope}
\end{tikzpicture}\hspace*{-0cm}
\par\end{center}

\begin{center}
(b)
\par\end{center}%
\end{minipage}
\par\end{centering}
\caption{(a) Two $E_{c}^{r}$'s and (b) the corresponding graphs}
\label{fig:twoErcs}
\end{figure}

Figure~\ref{fig:twoErcs}(a) shows two $E_{c}^{r}$'s with $k=6$.
So each of them have $6$ hexagons with centers (crosses) on another
hexagon.

Two circles in an $\erc$ may be disjoint, or intersect at one or
two points. To make out such relations between the circles in $\erc$,
a graph $G(\erc)=(\cV,\cE)$ can be used. Here the vertex set $\cV$
is simply $\Phi$ and the edge set is $\cE=\{(\lambda,\mu)\mid(\Phi r+\lambda c)\cap(\Phi r+\mu c)\not=\varnothing\}$.
For example, Figure~\ref{fig:twoErcs}(b) are the two graphs corresponding
to the two $E_{c}^{r}$'s on the left. In case that $(\lambda,\mu)\in\cE$,
the fact that $|(\Phi r+\lambda c)\cap(\Phi r+\mu c)|=1$ or $2$
is realized by coloring: an edge $(\lambda,\mu)\in\cE$ is \textit{even}
if $|(\Phi r+\lambda c)\cap(\Phi r+\mu c)|=2$ and \textit{odd} if
$|(\Phi r+\lambda c)\cap(\Phi r+\mu c)|=1$. For $j\in\{1,2,\dots,k-1\}$,
let $e_{j}=|(\Phi r+c)\cap(\Phi r+\varphi^{j}c)|$. Then the sequence
$e(r,c)=(\epsilon_{1},\epsilon_{2},\dots,\epsilon_{k-1})$ describes
completely the edge structure of $\erc$. (See \cite[(3.2)]{KeK96}.)

Abstractly, a sequence $e=(\epsilon_{1},\epsilon_{2},\dots,\epsilon_{k-1})$
with values $0$, $1$ and $2$ satisfying $\epsilon_{j}=\epsilon_{k-j}$
for $j=1,2,\dots,k-1$ gives rise to a colored graph $G(e)$. Here
$G(e)$ has the vertex set $\{v_{0},v_{2},\dots,v_{k-1}\}$ and edge
set
\[
\{(v_{i},v_{i+t})\mid0\leq i\leq k-1,\text{ }1\leq t\leq k/2,\text{ }\epsilon_{t}\not=0\}.
\]
An edge $(v_{i},v_{i+t})$ is even if $\epsilon_{t}=2$ and odd if
$\epsilon_{t}=1$. This way, one gets $G(\erc)=G(e(r,c))$. With such
abstraction, basic graphs can be defined. Let $j\in\bk=\{1,\dots,k-1\}$.
For the sequence $e=(\epsilon_{1},\epsilon_{2},\dots,\epsilon_{k-1})$
with $\epsilon_{i}=0$ if $i\not\in\{j,k-j\}$ and $\epsilon_{j}=\epsilon_{k-j}=1$,
the graph $\varGamma_{j}^{k}=G(e)$ is called the $j$th odd basic
$k$-graph. For the sequence $e=(\epsilon_{1},\epsilon_{2},\dots,\epsilon_{k-1})$
with $\epsilon_{i}=0$ if $i\not\in\{j,k-j\}$ and $\epsilon_{j}=\epsilon_{k-j}=2$,
the graph $\varPi_{j}^{k}=G(e)$ is called the $j$th even basic $k$-graph.
Specifically, if $\varDelta\in\{\varGamma_{j}^{k},\varPi_{j}^{k}\}$,
then the edge set $\cE(\varDelta)$ is $\{(v_{i},v_{i+j})\mid i\in\bk_{0}\}$
and $i+j$ is carried out modulo $k$.

\begin{figure}
\begin{centering}
\begin{minipage}[b][1\totalheight][t]{11.5cm}%
\hspace*{0.185cm}\begin{tikzpicture}[scale=0.85]
\newcommand\x{1.3cm}
\begin{scope}
\hhexagraphb
\end{scope}
\begin{scope}[xshift=2cm]
\node at (0,0) {$\longleftrightarrow$};
\end{scope}
\begin{scope}[xshift=4cm]
\hexabasica
\end{scope}
\begin{scope}[xshift=7cm]
\hexabasicb
\end{scope}
\begin{scope}[xshift=10cm]
\hexabasicc
\end{scope}
\end{tikzpicture}%
\end{minipage}
\par\end{centering}
\caption{An $G(E_{c}^{r})$ and the spanning basic graphs}
\label{fig:decomp_basic}
\end{figure}
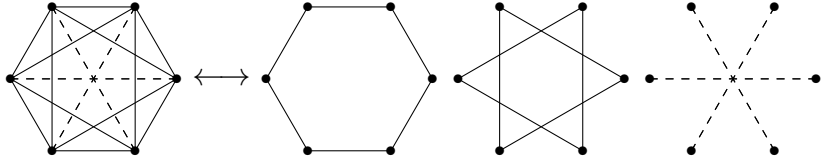

It turns out that each non-null graph $G(\erc)$ is the union of spanning
subgraphs, each of them is an even basic $k$-graph, or an odd basic
$k$-graph \cite[(4.1)]{KeK96}. Figure~\ref{fig:decomp_basic} shows
such a decomposition of the second graph in Figure~\ref{fig:twoErcs}(b)
into three basic graphs, two even ones (with solid-line edges) and
an odd one (with dotted-line edges).

Furthermore, an $i$th basic graph $\varDelta\in\{\varGamma_{i}^{k},\varPi_{i}^{k}\}$
is a spanning subgraph of $G(\erc)$ if and only if $c$ is in $\Phi(\varphi^{i}-1)^{-1}(\varphi^{j}-1)$
for some $j\in\bk$ \cite[(4.3)]{KeK96}. Thus, the set $\cM_{r}$
of $\erc$'s with $G(\erc)$ non-null is given by ${\cal M}_{r}=\{E_{c_{i,j}}^{r}\mid i,j\in\bk\}$,
where each $c_{i,j}=(\varphi^{i}-1)^{-1}(\varphi^{j}-1)$.

Finally, for $i\in\bk$, set $\gamma_{i}(r)=|\{\erc\in\cM_{r}\mid\varGamma_{i}^{k}\prec G(E_{c}^{r})\}|$
and $\pi_{i}(r)=|\{\erc\in\cM_{r}\mid\varPi_{i}^{k}\prec G(E_{c}^{r})\}|$,
where ``$\prec$'' means ``is a spanning subgraph of''. It was
shown in \cite[(4.7), (4.9), (4.10)]{KeK96} that
\begin{enumerate}
\item if $k$ is even, then $\gamma_{i}(r)=1$ and $\pi_{i}(r)=k/2-1$,
and\itemsep=-3pt
\item if $k$ is odd, then $\gamma_{i}(r)=k-1$ and $\pi_{i}(r)=0$.
\end{enumerate}
A natural question to ask now is what is the number of distinct graphs
in $\{G(\erc)\mid\erc\in\cM_{r}\}$? This amounts to learn when two
or more basic graphs are at the same time the spanning subgraphs of
some graph $G(\erc)$. In such a case, we say that \textit{these two
or more basic graphs overlap, }and that \textit{an overlap occurs
inside $G(\erc)$.}

From \cite[(4.3)]{KeK96}, one has
\begin{thm}
\label{thm:overlap_char}An overlap occurs inside $G(\erc)$ for some
$r,c\in F^{*}$ if and only if there exist $w\in\bk_{0}$, $i,j,s,t\in\bk$,
$i\not=s$, such that
\begin{equation}
(\varphi^{i}-1)^{-1}(\varphi^{j}-1)=\varphi^{w}(\varphi^{s}-1)^{-1}(\varphi^{t}-1).\label{E.olp}
\end{equation}
\end{thm}

\begin{rem}
\label{rem:trivial_overlaps}The situation $i=s$ and $j=k-t$ in
the theorem gives the same $i$th basic graph, and so actually no
overlap occurs. The situation $i=j$ and $s=t$ in the theorem describes
the overlap of the $i$th and the $s$th basic graphs in $G(E_{r}^{r})$.
Therefore, the graph $G(E_{r}^{r})$ is in fact a complete graph.
Thus, there is always an overlap with $\lfloor\frac{k}{2}\rfloor$
edges. In fact, this is the only overlap with $\lfloor\frac{k}{2}\rfloor$
edges. The number $\lfloor\frac{k}{2}\rfloor$ comes from the fact
that the edges $(v_{0},v_{i})$ and $(v_{0},v_{k-i})$ are the same
for all $i\in\bk$, as we have seen above. For obvious reasons, we'll
later refer to these overlaps as trivial.
\end{rem}

Our aim is to show in which situations overlaps can occur in $(F,\Phi)$.
The following two lemmas from \cite{KeK23-1} will be needed.
\begin{lem}[{{{{{\cite[Lemma 9]{KeK23-1}}}}}}]
\label{l:twocase} Let $(F,\Phi)$ be circular, and let $\chi=(\psi-1)^{-1}(\lambda-1)$
where $\lambda,\psi\in\Phi\setminus\{1\}$. If $\chi\in\Phi$, then
\[
\text{ either }\chi=1\text{ and }\psi=\lambda,\text{\quad or ~}\chi=-\lambda\text{ and }\psi=\lambda^{-1}.
\]
The second case implies either $p=2$, or that $|\Phi|$ is even.
\end{lem}

\begin{lem}[{{{{{\cite[Lemma 10]{KeK23-1}}}}}}]
\label{l:cjiPhi} Let $(p,k)$ be circular. For $i,j,t\in\bk$ we
have
\begin{enumerate}
\item if $k$ is even, then $c_{i,j}\in\Phi c_{i,t}\iff j=t\mbox{ or }i=k-t$,
and
\item if $k$ is odd, then $c_{i,j}\in\Phi c_{i,t}\iff j=t$ or, in case
$p=2$, $i=k-t$.
\end{enumerate}
\end{lem}

\section{Overlaps}

We shall fix $r\in F^{*}$ and for $c\in F^{*}$ denote $\varGamma_{c}=G(\erc)$.
Based on Theorem \ref{thm:overlap_char}, we make the following definition.
\begin{defn}
We say that the quadruple $(i,j\mid s,t)$, where $i,j,s,t\in\bk$,
$i\not=s$, forms an \textit{overlap} (with respect to $(F,\Phi)$,
or $(q,k)$) if $c_{i,j}\in\Phi c_{s,t}$. In this case, we also say
that $c_{i,j}$ is \textit{involved in an overlap}.
\end{defn}

We first collect some trivial cases, namely, an overlap $(i,j\mid s,t)$
with $j=i$, $j=t$ or $s=k-i$.

If $j=i$, then $c_{i,j}=1$, which puts $c_{s,t}$ into $\Phi=\Phi c_{j,i}$.
By Lemma~\ref{l:twocase}, if $t\not=s$, then either $2\div k$
and $t=k-s$, or $p=2$ and $t=k-s$.

If $j=t$, then
\[
\frac{\varphi^{j}-1}{\varphi^{i}-1}\in\Phi\frac{\varphi^{j}-1}{\varphi^{s}-1}\iff\frac{\varphi^{s}-1}{\varphi^{i}-1}=c_{i,s}\in\Phi=\Phi c_{j,j}
\]
so Lemma~\ref{l:twocase} applies again, and we have either $2\div k$
and $s=k-i$ where $i\not=\frac{k}{2}$, or $p=2$ and $s=k-i$ with
$i\not=s$.

If $s=k-i$, then we have $c_{i,j}\in\Phi c_{k-i,t}=\Phi c_{i,k-t}$
as well. Lemma~\ref{l:cjiPhi} says that if $t\not=k-j$, then either
$2\div k$ and $t=j$, or $p=2$ and $t=j$.

When describing or applying overlaps later, we mostly exclude the
instances above by referring to them as \textit{trivial overlaps}.
They are presented in compact form in Table~\ref{eq:trivOL}, where
the first column shows the forms of trivial overlaps and the second
column (if present) shows the extra conditions for the trivial overlaps
to occur.

\begin{table}[t]
\centering{}%
\begin{tabular}{c|c|c|c}
{\small$(i,i\mid s,s)$}  & {\small$(i,i\mid s,k-s)$}  & {\small$(i,j\mid k-i,j)$}  & {\small$(i,j\mid k-i,k-j)$}\tabularnewline
\hline
 & {\small$2\div k\vee p=2$}  & {\small$(2\div k\vee p=2)\wedge j\not=\frac{k}{2}$}  & {\small$i\not=\hk$}\tabularnewline
\end{tabular}\medskip{}
\caption[text]{{\small All trivial overlaps ($i,j,s\in\protect\bk$, $i\protect\not=s$)}}
\label{eq:trivOL}
\end{table}

Our aim is to determine the set of all nontrivial overlaps, namely,
\[
\cO=\cO(F,k)=\left\{ (i,j\mid s,t)\biggm|\Phi\frac{\vp^{j}-1}{\vp^{i}-1}=\Phi\frac{\vp^{t}-1}{\vp^{s}-1},\text{}(i,j,s,t)\in\cI\right\} .
\]
As $s=t$ would create a trivial overlap (see Lemma~\ref{l:cjiPhi}),
we have added the condition $s\neq t$ to the above definition of
$\cO$ for symmetry. Also, we could have used the notation $\cO(F,\Phi)$
instead of $\cO(F,k)$, but $\Phi$ is uniquely determined by $k$.
As we will also consider the complex number field $\bC$, this notation
comes in handy. Actually, the biggest part of this chapter will be
occupied by the case when $F=\bC$.

After clearing denominators and expanding \eqref{E.olp}, we find
\begin{equation}
\vp^{\omega+j+s}+\vp^{\omega}-\vp^{\omega+s}-\vp^{\omega+j}-\vp^{t+i}+\vp^{t}+\vp^{i}-1=0,\mbox{~with \ensuremath{i\not=s}}.\label{e:poly1}
\end{equation}
For $i,j,s,t\in\bk$, $i\not=s$, and $\o\in\bk_{0}$ define polynomials
over $F$
\begin{align}
f_{i,j,s,t,\o}(x) & =x^{\omega+j+s}+x^{\omega}-x^{\omega+s}-x^{\omega+j}-x^{t+i}+x^{t}+x^{i}-1\label{e:poly}\\
 & =x^{\omega}(x^{j}-1)(x^{s}-1)-(x^{i}-1)(x^{t}-1).\label{e:poly-circ}
\end{align}
Obviously, when $i\not=s$, then we have $(i,j\mid s,t)$ is an overlap
if and only if there exists some $\o\in\bk_{0}$ such that $f_{i,j,s,t,\o}(\varphi)=0$.

A direct consequence of this is
\begin{lem}
\label{l:Fextend} Let $(F,\Phi)$ be circular, and let $K$ be an
extension field of $F$. Then $(K,k)$ is circular and $\cO(K,k)=\cO(F,k)$.
\end{lem}

\begin{proof}
The statement about circularity comes directly from Theorem~\ref{l:circ}.
The second statement is clear.
\end{proof}
This means that we can reduce our discussions to the smallest subfield
of $F$ containing a $k$-th root of unity. In particular, if $F$
is finite, the set $\cO(F,k)$ only depends on the characteristic~$p$.
We therefore sometimes simply write $\cO(p,k)$ for $\cO(F,k)$.

\subsection{The Complex Numbers}

\label{sec:complex-numbers}

We now set the stage for $\bC$. Let $\phi=\exp(2\pi\bm{i}/k)$ be
a primitive $k$-th root of unity in $\bC$ where $\bm{i}^{2}=-1$,
and let $U=U_{k}=\gen[\phi]$. Notice that $(\bC,U)$ is circular
for all $k$ as $U$ is a subset of the unit circle. We will prove
\begin{thm}
\label{th:OLPrimes} Let $p$ be a prime. Then $\cO(\bC,k)\subseteq\cO(p,k)$.
Moreover, for $k\geq3$, the set $\cQk=\{p\mbox{ prime}\mid p\text{ divides }k\text{ or }\cO(\bC,k)\neq\cO(p,k)\}$
is finite.
\end{thm}

To prepare the proof we first show a proposition. All information
on cyclotomic fields needed for this can be found in \cite[Ch.\ 13, \S2]{IreRos}.
\begin{prop}
\label{l:polyCyc} Let $p$ be a prime, $k\in\N$, $k\ge3$, and $F$
a field of characteristic~$p$ which contains an element of order
$k$ in $F^{*}$. Let $\phi$ be a primitive $k$-th root of unity
over $\bQ$ and let $\cN:\bQ(\phi)\to\bQ$, be the Galois norm. Let
$g$ be a polynomial in $\bZ[x]$.
\begin{enumerate}
\item There exists an element $\varphi$ of order $k$ in $F^{*}$ such
that $g(\varphi)=0$ if and only if $p\div\cN(g(\phi))$.\itemsep=-3pt
\item In case $g(\phi)=0$ (which is equivalent to $\cN(g(\phi))=0$), we
have $g(\psi)=0$ for all elements $\psi$ of order $k$ in $F^{*}$.
\end{enumerate}
\end{prop}

\begin{proof}
There is no loss in generality to assume that $F$ is the smallest
field with the given properties. As $\bZ[\phi]$ is the ring of integers
inside the $k$th cyclotomic field $\bQ(\phi)$, there exists a ring-epimorphism
$\theta:\bZ[\phi]\to F;u\mapsto u^{\theta}$ mapping $\phi$ to some
primitive $k$-th roots of unity, i.e., an element of order $k$,
inside $F$. (If $F$ where not smallest, the image of $\theta$ would
be this smallest field.) We note that the kernel $P$ of $\theta$
is a prime ideal containing $p$, and the map $\theta$ extends naturally
to a ring-epimorphism of the polynomial rings $\bZ[\phi][x]\to F[x]$,
which is also denoted by $\theta$. Let $g\in\bZ[x]$ be a polynomial
and let $G$ be the Galois group of $[\bQ(\phi){:}\bQ]$. As norms
of elements of $\bZ[\phi]$ are all in $\bZ$, we have that $\cN(g(\phi))\in\bZ$.
Thus we find
\begin{align*}
p\div\cN(g(\phi)) & \iff\cN(g(\phi))=\prod_{\sigma\in G}(g(\phi))^{\sigma}=\prod_{\sigma\in G}g\left(\phi^{\sigma}\right)\in P\\
 & \iff g\left(\phi^{\sigma_{0}}\right)\in P\text{\quad for some \ensuremath{\sigma_{0}\in G}}\\
 & \iff\left(g\left(\phi^{d}\right)\right)^{\theta}=0\text{\quad for some \ensuremath{d\in\bZ_{k}^{\times}}}
\end{align*}
since $P$ is a prime ideal, and the group $G$ is naturally isomorphic
to the group of units $\bZ_{k}^{\times}$ of the ring $\bZ_{k}$.
(Indeed, if $\sigma_{0}$ and $d$ correspond under this isomorphism,
then $\phi^{\sigma_{0}}=\phi^{d}$.)

Now, assume that $p\div\cN(g(\phi))$. Then there exists a $d\in\bZ_{k}^{\times}$
such that $\varphi=\left(\phi^{d}\right)^{\theta}$ is a root of $g$.
Conversely, assume that $\varphi$ exists. Then there exists a preimage
of $\varphi$ under $\theta$, which is a primitive $k$-th roots
of unity in $\bQ(\phi)$. As $G$ is transitive on the primitive $k$-th
roots of unity, there exists $d\in\bZ^{\times}$ with $\varphi=\left(\phi^{d}\right)^{\theta}.$
Therefore, $\left(g\left(\phi^{d}\right)\right)^{\theta}=g(\varphi)=0,$
and so $p\div\cN(g(\phi))$. This proves (1).

Next, suppose that $g(\phi)=0$. We have $g(\phi^{d})=0$ for all
$d\in\bZ_{k}^{\times}$. Therefore,
\[
0=g(\phi^{d})^{\theta}=g((\phi^{\theta})^{d}).
\]
Here, $(\phi^{\theta})^{d}$, $d\in\bZ_{k}^{\times}$, are exactly
the elements of order $k$ inside $F^{*}$. This is (2).
\end{proof}
\begin{rem}
\label{rem:polyCyc} In the above proof, $\cN(g(\phi))=0$ implies
that there exists a conjugate $\phi^{d_{0}}$ of $\phi$ such that
$g(\phi^{d_{0}})=0$. But then, by the action of $G$, $g(\phi^{d})=0$
for all $d\in\bZ^{\times}$, and thus $g(\varphi)=0$ for all elements
$\varphi\in F$ of order $k$.
\end{rem}

\begin{proof}[Proof of Theorem~\ref{th:OLPrimes}]
We shall use the Galois norm $\cN$ as introduced in Proposition~\ref{l:polyCyc}.

For $(i,j\mid s,t)\in\cO(\bC,k)$ there exists $\o\in\bk_{0}$ such
that $f_{i,j,s,t,\o}(\phi)=0$. Then $(i,j,s,t)\in\cI$ and $\cN(f_{i,j,s,t,\o}(\phi))=0$,
thus $(i,j\mid s,t)\in\cO(p,k)$ by Proposition~\ref{l:polyCyc}.

Suppose that $p\in\cQ_{k}$ and $p\nmid k$. Thus there exists $(i,j\mid s,t)\in\cO(p,k)\setminus\cO(\bC,k)$
with a corresponding $\o$. By Proposition~\ref{l:polyCyc} again,
$\cN(f_{i,j,s,t,\o}(\phi))\not=0$ and $p\div\cN(f_{i,j,s,t,\o}(\phi))$.
There are only finitely many polynomials~$f_{i,j,s,t,\o}$ and each
integer $\cN(f_{i,j,s,t,\o}(\phi))$ has only finitely many prime
divisors. Thus the set $\cQk$ is finite as well.
\end{proof}
\begin{rem}
$f_{i,j,s,t,\o}(\phi)$ is a sum of $8$ roots of unity. Let $e$
be the number of elements in $\bk$ coprime to $k$. To form $\cN(g(\phi))$,
we multiply $e$ such sums. After expansion, we have a total of~$8^{e}$
summands, each of which is a product of roots of unity and has absolute
value $1$. This yields the inequality $\abs[\cN(g(\phi))]\le8^{e}$.
Thus for every prime $p\in\cQk$ we have $p<8^{e}$. Note that this
is a very crude bound as our data show, and is suggested from the
proof, too.
\end{rem}

For every $k\ge3$, the set $\cQk$ is referred to as the set of \textit{exceptional
primes}.
\begin{cor}
\label{PksubsetQk}It holds that $\cPk\subseteq\cQk$. Consequently,
if $(p,k)$ is a Ferrero pair with $p\not\in\cQk$, then $(p,k)$
is circular.
\end{cor}

\begin{proof}
We notice that the polynomials (\ref{eq:fijst}) used for getting
$\cPk$ are all of the form $f_{i,j,s,t,0}(x)$, $(i,j,s,t)\in\cI$.
As $(\mathbb{C},U)$ is circular, $f_{i,j,s,t,0}(\phi)\not=0$ for
all $(i,j,s,t)\in\cI$ by Lemma~\ref{l:circ}. From Proposition~\ref{l:polyCyc},
it follows immediately that $\cPk$ consists exactly the prime divisors
of ${\cal N}(f_{i,j,s,t,0}(\phi))$, $(i,j,s,t)\in\cI$. Thus, in
both $\cQk$ and $\cPk$, we are determining primes dividing $\cN(f_{i,j,s,t,\o}(\phi))$
wherever $\cN(f_{i,j,s,t,\o}(\phi))\not=0$ for $(i,j,s,t)\in\cI$:
in the case of $\cQk$, $\o\in\bk_{0}$ and in the case of $\cPk$,
$\o=0$. Therefore, we have $\cPk\subseteq\cQk$.
\end{proof}
We provide some examples of $\cQk$ with elements of $\cPk$ underlined
in Table~\ref{tab:Qk}. The algorithm to find the elements is based
on the above proof.
\begin{table}
\begin{centering}
\begin{tabular}{cccl}
\hspace*{-1cm}  & {\small$\cQ_{4}$}  & \hspace*{-3mm}$=$  & \hspace*{-2.5mm}{\small$\{\underline{2},3,\underline{5}\},$}\tabularnewline
 & {\small$\cQ_{5}$}  & \hspace*{-3mm}$=$  & \hspace*{-2.5mm}{\small$\{\underline{5},\underline{11}\},$}\tabularnewline
 & {\small$\cQ_{6}$}  & \hspace*{-3mm}$=$  & \hspace*{-2.5mm}{\small$\{\underline{2},\underline{3},5,\underline{7},\underline{13},\underline{19},31,37\},$}\tabularnewline
 & {\small$\cQ_{7}$}  & \hspace*{-3mm}$=$  & \hspace*{-2.5mm}{\small$\{\underline{2},\underline{7},13,\underline{29},\underline{43},71\},$}\tabularnewline
 & {\small$\cQ_{8}$}  & \hspace*{-3mm}$=$  & \hspace*{-2.5mm}{\small$\{\underline{2},\underline{3},\underline{5},{7},13,\underline{17},\underline{41},73,89,97,113\},$}\tabularnewline
 & {\small$\cQ_{9}$}  & \hspace*{-3mm}$=$  & \hspace*{-2.5mm}{\small$\{2,\underline{3},17,\underline{19},\underline{37},\underline{73},\underline{109},\underline{127},163,181,199,\underline{271},397,541\},$}\tabularnewline
 & {\small$\cQ_{10}$}  & \hspace*{-3mm}$=$  & \hspace*{-2.5mm}{\small{}{}{}{}{}}{\small{}%
\begin{minipage}[t]{0.8\columnwidth}%
{\small$\{\underline{2},3,\underline{5},\underline{11},19,29,\underline{31},\underline{41},\underline{61},\underline{71},\underline{101},131,151,181,191,211,241,$}\\
 {\small$\text{\enspace}$$251,271,281,311,331,401,421,541,641,761,881,941\},$}%
\end{minipage}}\tabularnewline
 & {\small$\cQ_{11}$}  & \hspace*{-3mm}$=$  & \hspace*{-2.5mm}{\small{}{}{}{}{}}{\small{}%
\begin{minipage}[t]{0.8\columnwidth}%
{\small$\{3,\underline{11},\underline{23},43,\underline{67},\underline{89},109,\underline{199},331,\underline{353},\underline{397},419,463,617,661,$}\\
 {\small$\text{\enspace}$$\underline{683},727,859,881,947,991,1277,1453,2069,2311,2399\},$}%
\end{minipage}}\tabularnewline
 & {\small$\cQ_{12}$}  & \hspace*{-3mm}$=$  & \hspace*{-2.5mm}{\small{}{}{}{}{}}{\small{}%
\begin{minipage}[t]{0.8\columnwidth}%
{\small$\{\underline{2},\underline{3},\underline{5},\underline{7},11,\underline{13},\underline{17},\underline{19},23,29,31,\underline{37},\underline{61},\underline{73},\underline{97},\underline{109},\underline{157},\underline{181},$}\\
 {\small$\text{\enspace}$$\underline{193},229,241,277,313,337,349,373,397,409,421,433,541,601,$}\\
 {\small$\text{\enspace}$$661,769,1009\}.$}%
\end{minipage}}\tabularnewline
\end{tabular}\medskip{}
\par\end{centering}
\caption[]{{\small Exceptional primes; elements from $\protect\cPk$ underlined}}
\label{tab:Qk}
\end{table}

The following examples from \cite{KeK96} are some nontrivial overlaps
for Ferrero pairs $(q,k)$. Note that these are universal in the sense
that they do not depend on $q$ (or $p$), but only on the shape of~$k$.
\begin{example}
\label{ex:overlap} If $k=6\ell$, $\ell\in\N$, then $\varphi^{\ell}$
is a sixth root of unity, and $\varphi^{3\ell}=-1=\varphi^{2\ell}-\varphi^{\ell}$.
Therefore
\begin{align*}
\varphi^{\ell}\frac{\varphi^{\ell}-1}{\varphi^{i}-1} & =\frac{-1}{\varphi^{i}-1}\\
\noalign{\hbox{and}}\varphi^{i}\frac{\varphi^{3\ell-i}-1}{\varphi^{2i}-1} & =\frac{\varphi^{3\ell}-\varphi^{i}}{(\varphi^{i}-1)(\varphi^{i}+1)}=\frac{-1-\varphi^{i}}{(\varphi^{i}-1)(\varphi^{i}+1)}=\frac{-1}{\varphi^{i}-1}.
\end{align*}
This yields $c_{2i,3\ell-i}=\varphi^{\ell-i}c_{i,\ell}$ for all $1\le i\leq k/4$.
To put it short, we have that $(i,\ell\mid2i,3\ell-i)$ forms an overlap
for every $i\in\{1,\dots,\lfloor k/4\rfloor\}$.

Notice that the case $i=\ell$ is trivial, but all other cases are
not. Thus nontrivial examples of this kind start with $k=12$.
\end{example}

As we represent overlaps by the exponents with respect to a fixed
generator, the actual quadruples will depend on this generator. We
give examples for this in Examples~\ref{ex:p13k7} and~\ref{ex:p11k12}.

The main concern of this paper is the determination of $\cO(\bC,k)$
and thus by Theorem~\ref{th:OLPrimes} that of $\cO(F,k)$ for all
finite fields with characteristic not in $\cQk$. We will now show
that the set does not really depend on the generator in this case.
In other words: a problem occurs only for exceptional primes.
\begin{lem}
\label{l:uvstIndep} Let $k\ge3$ and let $F$ be a field of characteristic
$p\notin\cQ_{k}$ such that $F^{*}$ contains a subgroup~$\Phi$
of order $k$, or $F=\bC$. Then the set $\cO(F,k)$ is independent
of the choice of the generator for $\Phi$. Specifically, let $\psi$
and $\chi$ be generators of $\Phi$, and let $(i,j\mid s,t)\in\cO(F,k)$,
i.e. there exists $w\in\bk_{0}$ such that
\[
\psi^{w}\cdot\frac{\psi^{j}-1}{\psi^{i}-1}=\frac{\psi^{s}-1}{\psi^{t}-1}\Longleftrightarrow\chi^{w}\cdot\frac{\chi^{j}-1}{\chi^{i}-1}=\frac{\chi^{s}-1}{\chi^{t}-1}.
\]
\end{lem}

\begin{proof}
We first treat the complex case. We can restrict to the $k$th cyclotomic
field $F=\bQ(\psi)$. There exists an automorphism $\sigma$ of $F$
such that $\psi^{\sigma}=\chi$. The first equation implies $f_{i,j,s,t,\o}(\psi)=0$,
then also
\[
0=f_{i,j,s,t,\o}(\psi)^{\sigma}=f_{i,j,s,t,\o}\left(\psi^{\sigma}\right)=f_{i,j,s,t,\o}(\chi).
\]
Now, the finite case follows directly with Proposition~\ref{l:polyCyc}.
\end{proof}

\subsection{The Reduced Form}

\label{sec:reduced-form}

To reduce complexity of the set $\cO(F,k)$ we use some group actions
on $\cO(F,k)$. We consider the mappings $\kappa_{\ell}:\bk^{4}\to\bk^{4}$,
$\ell\in\{1,2,3,4\}$, which transforms the $\ell$-th entry $u_{\ell}$
of a quadruple to $k-u_{\ell}$. These four mappings generate an elementary
abelian $2$-group $\cK_{0}$ of order $16$. The subgroup $\cK_{1}=\gen[\kappa_{1}\kappa_{4},\kappa_{2}\kappa_{4},\kappa_{3}\kappa_{4}]$
generated by products of two such generators has index $2$ in~$\cK_{0}$.

Let $(F,k)$ be a circular Ferrero pair, then $\cK_{1}$ acts on $\cO(F,k)$.
If $k$ is even, then $\cK_{0}$ acts on $\cO(F,k)$ since $-1=\varphi^{\hk}\in\Phi$.

If $o:=(i,j\mid s,t)\in\cO(F,k)$ is an overlap then the following
permutations of the entries of $o$ give more identities as in (\ref{E.olp})
\begin{equation}
(i,t),\ (j,s),\ (i,t)(j,s),\ (i,s)(j,t),\ (i,j)(s,t),\ (i,j,t,s),\ (i,s,t,j).\label{e:perms}
\end{equation}
These together with the identity map form a dihedral group $D_{4}$
acting on $\cO(F,k)$, too. It is easy to see that $D_{4}$ normalizes
$\cK_{1}$ (and also $\cK_{0}$). Thus the semidirect product of $D_{4}$
together with $\cK_{1}$, or $\cK_{0}$ form groups $\cG_{1}$, or
$\cG_{0}$, respectively.
\begin{lem}
\label{l:Gaction} Let $(F,k)$ be a circular Ferrero pair. If $k$
is odd, then $\cG_{1}$ acts on $\cO(F,k)$, and if $k$ is even,
then $\cG_{0}$ acts on $\cO(F,k)$.
\end{lem}

Occasionally, we will write $o\sim o'$ if two overlaps $o,o'\in\cO(F,k)$
are related by the group action of Lemma~\ref{l:Gaction}. Clearly,
$\sim$ is an equivalence relation on the set of all overlaps. To
describe this set it suffices to give a representative for each class.
We will now describe a ``reduced'' representative for each class.
Whenever situation allows, we will choose a reduced representative
in our exposition.

Let $(i,j\mid s,t)\in\cO(F,k)$ be a nontrivial overlap. By applying
elements from $\cK_{1}$, we can pass to an equivalent quadruple which
has at most one entry greater than $\hk$. If $k$ is even, we can
even pass to an equivalent quadruple which has all of its entries
less than or equal to~$\hk$.

By applying the permutations $(i,t)$ if necessary, we can assume
that $i\le t$. Applying $(j,s)$ and/or $(i,j)(s,t)$ we may assume
that $i$ is the smallest of the values among $i,j,s,t$. Now, applying
$(j,s)$, again, if necessary, we can assume with no loss of generality
that
\begin{equation}
i<j\leq s\mbox{ \ and \ }j\le\hk.\label{eq:uvs}
\end{equation}
An element $(i,j\mid s,t)\in\cO(F,k)$ is called \textit{reduced}
if it satisfies the conditions in \eqref{eq:uvs} and has at most
one entry greater than $\hk$.
\begin{lem}
In each equivalence class of $\cO(F,k)$ there exists a reduced element.
\end{lem}

\begin{rem}
The reduced form is not unique. E.g., let $k=12$. From Example~\ref{ex:overlap}
we have an overlap $o=(1,2\mid2,5)$, which clearly is reduced. However,
$(1,2\mid2,7)\sim o$ is also reduced; and so is $(1,2\mid10,5)\sim o$.
\end{rem}

We emphasis again that for a prime $p\in\cQ_{k}$, the representation
of the overlaps as powers in the set $\cO(F,k)$ depends on the choice
of the generator for~$\Phi$.
\begin{example}
\label{ex:p13k7} For $p=13$ and $k=7$, we have $13\in\cQ_{7}\setminus\cP_{7}$.
Thus $(13^{2},7)$ is circular, but there exist exceptional overlaps
such as $(1,2\mid2,6)\sim(1,2\mid5,1)$ both of which are reduced.
This overlap works with the element $\varphi$ of order $k$ in the
quadratic extension of $\bZ_{p}$ with minimal polynomial $x^{2}+3x+1$.\footnote{The minimal polynomial of the other generators are $x^{2}+5x+1$ (for
$\vp^{\pm3}$) and $x^{2}+6x+1$ (for $\vp^{\pm2}$).} Indeed, a simple computation shows that
\begin{alignat*}{1}
\vp^{5}\frac{\vp^{2}-1}{\varphi-1}=\frac{\vp^{6}-1}{\vp^{2}-1} & \iff(\vp^{2}-1)^{2}\vp^{-2}=(\vp^{-1}-1)(\varphi-1)\\
 & \iff\vp^{2}+\vp^{-2}-2=2-(\vp^{-1}+\varphi).
\end{alignat*}
The last equation holds since the trace of $\varphi$ is $-3$ and
that of $\vp^{2}$ is $-6$.

Another overlap which works with $\varphi$ is $(1,2\mid5,1)$. In
the same way we have $(1,3\mid3,6)\sim(1,3\mid4,1)$ working with~$\vp^{2}$.
Notice that, however, $(1,3\mid3,6)$ does not work with $\varphi$.
\end{example}

\begin{example}
\label{ex:p11k12} For $p=11$ and $k=12$ we again have $11\in\cQ_{12}\setminus\cP_{12}$.
Thus $(11^{2},12)$ is circular. Besides the natural overlaps $(1,2\mid2,5)$
and $(2,3\mid3,6)$ from Theorem~\ref{th.main} there exist exceptional
overlaps such as $(2,4\mid5,3)$ which works with $\varphi$, a root
of $x^{2}+5x+1$ (and $\o=4$). On the other hand, $(1,2\mid3,4)$
works with $\vp^{5}$, which has minimal polynomial $x^{2}-5x+1$.
\end{example}

\section{Normalized Form}

\label{S:normalized}

We will specialize to the realm of the complex numbers and consider
only nontrivial overlaps. It turns out, as we shall see later, that
there can be only trivial overlaps when $k$ is odd. \textbf{From
now on, we shall assume that $k$ is even.} We come back to the odd
case only in Theorem~\ref{th:OLodd}. By abuse of notation, we will
write $\cO=\cO(\bC,k)$.

As the group $\cG_{0}$ acts on $\cO$ for even $k$, we can assume
that $i,j,s,t\in\{1,2,\dots,\hk\}$. Recall that $\phi=\exp(2\pi\bm{i}/k)$
is a primitive $k$-th root of unity in $\bC$. We will use the polar
decomposition of $\vpc[r]-1$, $r\in\mathbb{R}$. This is easily computed
using the identity $\vpc[r]-1=\left(\vpc[\frac{r}{2}]-\vpc[-\frac{r}{2}]\right)\vpc[\frac{r}{2}]$
and Euler's formula.
\begin{lem}
\label{l:polardec}For $r\in\mathbb{R}$, it holds that
\[
{\displaystyle \vpc[r]-1=2\sin\frac{\pi r}{k}\cdot\exp\bm{i}\left(\frac{\pi}{2}+\frac{\pi r}{k}\right)}.
\]
\end{lem}

The following observations further reduces the overlap quadruples
of interest.
\begin{lem}
\label{l:order-uv}Let $(i,j\mid s,t)\in\cO$ with $i,j,s,t\in\{1,2,\dots,\frac{k}{2}\}$.
Then
\[
i<j\iff s<t\text{ and }i<s\iff j<t.
\]
\end{lem}

\begin{proof}
By Lemma~\ref{l:polardec} and the monotonicity of sine on the interval
$[0,\frac{\pi}{2}]$, we have
\[
\abs[\frac{\vpc[j]-1}{\vpc[i]-1}]=\frac{\sin\frac{\pi j}{k}}{\sin\frac{\pi i}{k}}>1\iff i<j.
\]
As the same holds for $(s,t)$, the first statement of the lemma follows.
By exchanging the roles of $j$ and $s$, the second statement follows
immediately.
\end{proof}
Therefore (\ref{eq:uvs}) implies $s<t$ for a reduced quadruple $(i,j\mid s,t)$,
$i,j,s,t\in\{1,2,\dots,\frac{k}{2}\}$. Summarizing we can assume
\begin{equation}
0<i<j\leq s<t\leq\hk.\label{eq:normal-O}
\end{equation}
An element $(i,j\mid s,t)\in\cO$ is called \textit{normalized} if
it satisfies the conditions in \eqref{eq:normal-O}.

With these premises we find
\begin{lem}
\label{l:sumless} If $(i,j\mid s,t)\in\cO$ is normalized, then $j-i<t-s$,
and so $j+s<i+t<k$.
\end{lem}

\begin{proof}
It suffices to compare the absolute values of the left and right hand
side of (\ref{E.olp}). This gives, by Lemma~\ref{l:polardec},
\[
\frac{\sin\frac{\pi j}{k}}{\sin\frac{\pi i}{k}}=\frac{\sin\frac{\pi t}{k}}{\sin\frac{\pi s}{k}}.
\]
Set $f(x,y)=\frac{\sin(x+y)}{\sin x}$ on the set $\{(x,y)\mid0<x<x+y\le\frac{\pi}{2}\}$.
Simple calculus analysis reveals that, keeping $y$ fixed, $f(x,y)$
strictly decreases as $x$ increases, and, keeping $x$ fixed, it
strictly increases as $y$ increases.

Now, from
\begin{multline*}
\qquad f\left(\frac{\pi}{k}s,\frac{\pi}{k}(j-i)\right)<f\left(\frac{\pi}{k}i,\frac{\pi}{k}(j-i)\right)\\
=\frac{\sin\frac{\pi j}{k}}{\sin\frac{\pi i}{k}}=\frac{\sin\frac{\pi t}{k}}{\sin\frac{\pi s}{k}}=f\left(\frac{\pi}{k}s,\frac{\pi}{k}(t-s)\right)\qquad
\end{multline*}
we infer that $j-i<t-s$.
\end{proof}
\begin{rem}
A reduced quadruple is not necessarily normalized, but the converse
is true. Over finite fields there do exist reduced quadruples which
cannot be normalized as Example~\ref{ex:p11k12} shows. This phenomenon
can only occurs if the characteristic of the field is an exceptional
prime.
\end{rem}

\subsection{Beyond the Normalized Form}

Not all permutations in $D_{4}$ viewed as a subgroup of $\cG_{0}$
give distinct elements from $\cO$. Indeed
\begin{lem}
If $o:=(i,j\mid s,t)\in\cO$ is normalized then there are at most
four values modulo $\Phi$ derived from this by permutations, namely,
\[
\frac{\vpc[j]-1}{\vpc[i]-1},~\frac{\vpc[s]-1}{\vpc[i]-1},~\frac{\vpc[i]-1}{\vpc[j]-1}\text{ and }\frac{\vpc[i]-1}{\vpc[s]-1}.
\]
The action of $\cK_{0}$ does not change the cosets.
\end{lem}

\begin{proof}
Starting from the first value in the theorem, the first permutation
from \eqref{e:perms} produces the last entry, the second and third
produce the second and third entry, respectively.

The permutation $(i,s)(j,t)$ only interchanges the left hand side
with the right hand side of \eqref{E.olp} up to an factor in $\Phi$.
Therefore the other three permutations in \eqref{e:perms} cannot
give more solutions either.
\end{proof}

\section{Main Theorem}

\label{sec:Main-Theorem}

Finally, in this section we reach our principal goal, the determination
of the set $\cO=\cO(\mathbb{C},k)$ of nontrivial overlaps over $\bC$.
For any finite field $F$ with characteristic not in $\cQ_{k}$ this
set coincides with $\cO(F,k)$, cf.\,Theorem~\ref{th:OLPrimes}.

Later in Theorem~\ref{t:tripolp} we also determine the triple overlaps
and prove that there are no quadruple overlaps. For the sake of easy
reference, we include the findings of the triple overlaps in Theorem~\ref{t:tripolp}
into the following theorem (the last column, marked with~$T_{r}$).
\begin{thm}
\label{th.main} Let $k$ be even and $\cO$ nonempty, then there
exists $\ell\in\N$ such that $k=6\ell$. Depending on the shape of
$k$, $\cO$ is a union of the corresponding sets $\cO_{1}$, $\cO_{30}$,
$\cO_{42}$, and $\cO_{60}$ described below. When we write $k=N\ell_{r}$
for $N\in\{30,42,60\}$, we mean that $k$ is divisible by $N$.\LetLtxMacro\itemold\item
\renewcommand{\item}{\itemindent-0.8cm\itemold}
\begin{enumerate}
\item For $k=6\ell$, we have
\[
\p{\ell-u}{\ell}{u}{3\ell-u}{2u}
\]
{\small${}$\hspace{-0.8cm}}for $1\le u\leq\big\lfloor\tfrac{k}{4}\big\rfloor$
with $u\not=\frac{k}{6}=\ell$. Namely, $\cO_{1}$ consists of those\linebreak{}
 {\small${}$\hspace{-0.8cm}}$(i,j,s,t)\sim(u,\ell\mid2u,3\ell-u)$,
where $1\le u\leq\bigl\lfloor\tfrac{k}{4}\bigr\rfloor$ and $u\not=\tfrac{k}{6}=\ell$.
\\
\item For $k=30\ell_{1}$, the normalized forms of the elements in $\cO_{30}$
are \medskip{}
 \\
{\small${}$\hspace{-1cm}}{\small{}%
\begin{tabular}{lcl}
{\small$\plist[3\ell_{1}]{3\ell_{1}}{3\ell_{1}}{\ell_{1}}{11\ell_{1}}$,
 } &  & \tabularnewline
{\small$\plist[5\ell_{1}]{\ell_{1}}{5\ell_{1}}{3\ell_{1}}{9\ell_{1}}$,
 } &  & {\small\kern-.6cm $T_{1}:(3\ell_{1},5\ell_{1}\mid5\ell_{1},9\ell_{1}\mid6\ell_{1},12\ell_{1})$,}\tabularnewline
{\small$\plist[9\ell_{1}]{\ell_{1}}{9\ell_{1}}{7\ell_{1}}{13\ell_{1}}$,
 } &  & \tabularnewline
 &  & \tabularnewline
{\small$\plist[4\ell_{1}]{2\ell_{1}}{2\ell_{1}}{\ell_{1}}{9\ell_{1}}$,
 } &  & {\small\kern-.6cm $T_{2}:(\ell_{1},2\ell_{1}\mid4\ell_{1},9\ell_{1}\mid5\ell_{1},14\ell_{1})$,}\tabularnewline
{\small$\plist[5\ell_{1}]{\ell_{1}}{3\ell_{1}}{2\ell_{1}}{8\ell_{1}}$,
 } &  & {\small\kern-.6cm $\ensuremath{T_{3}}:\ensuremath{(2\ell_{1},3\ell_{1}\mid5\ell_{1},8\ell_{1}\mid7\ell_{1},14\ell_{1})}$,}\tabularnewline
 &  & {\small\kern-.6cm $\ensuremath{T_{4}}:\ensuremath{(2\ell_{1},5\ell_{1}\mid3\ell_{1},8\ell_{1}\mid4\ell_{1},13\ell_{1})}$,}\tabularnewline
{\small$\plist[11\ell_{1}]{\ell_{1}}{9\ell_{1}}{8\ell_{1}}{14\ell_{1}}$,
 } &  & {\small\kern-.6cm $\ensuremath{T_{5}}:\ensuremath{(4\ell_{1},5\ell_{1}\mid8\ell_{1},11\ell_{1}\mid9\ell_{1},14\ell_{1})}$,}\tabularnewline
{\small$\plist[7\ell_{1}]{3\ell_{1}}{3\ell_{1}}{2\ell_{1}}{14\ell_{1}}$,
 } &  & {\small\kern-.6cm $\ensuremath{T_{3}}:\ensuremath{(2\ell_{1},3\ell_{1}\mid5\ell_{1},8\ell_{1}\mid7\ell_{1},14\ell_{1})}$,}\tabularnewline
{\small$\plist[8\ell_{1}]{2\ell_{1}}{4\ell_{1}}{3\ell_{1}}{13\ell_{1}}$,
 } &  & {\small\kern-.6cm $\ensuremath{T_{4}}:\ensuremath{(2\ell_{1},5\ell_{1}\mid3\ell_{1},8\ell_{1}\mid4\ell_{1},13\ell_{1})}$,}\tabularnewline
{\small$\plist[9\ell_{1}]{2\ell_{1}}{5\ell_{1}}{4\ell_{1}}{14\ell_{1}}$,
 } &  & {\small\kern-.6cm $\ensuremath{T_{2}}:\ensuremath{(\ell_{1},2\ell_{1}\mid4\ell_{1},9\ell_{1}\mid5\ell_{1},14\ell_{1})}$,}\tabularnewline
 &  & {\small\kern-.6cm $\ensuremath{T_{5}}:\ensuremath{(4\ell_{1},5\ell_{1}\mid8\ell_{1},11\ell_{1}\mid9\ell_{1},14\ell_{1})}$.}\tabularnewline
\end{tabular}}{\small}\\
{\small\par}
\item For $k=42\ell_{2}$, the normalized forms of the elements in $\cO_{42}$
are \medskip{}
 \\
{\small${}$\hspace{-1cm}}{\small{}%
\begin{tabular}{lcc}
{\small$\plist[9\ell_{2}]{3\ell_{2}}{3\ell_{2}}{2\ell_{2}}{16\ell_{2}}$,
 } &  & \tabularnewline
{\small$\plist[10\ell_{2}]{2\ell_{2}}{4\ell_{2}}{3\ell_{2}}{15\ell_{2}}$,
 } &  & \tabularnewline
{\small$\plist[15\ell_{2}]{2\ell_{2}}{9\ell_{2}}{8\ell_{2}}{20\ell_{2}}$.
 } &  & \tabularnewline
\end{tabular}}{\small}\\
{\small\par}
\item For $k=60\ell_{3}$, the normalized forms of the elements in $\cO_{60}$
are \medskip{}
 \\
{\small${}$\hspace{-1cm}}{\small{}%
\begin{tabular}{lcc}
{\small$\plist[16\ell_{3}]{5\ell_{3}}{4\ell_{3}}{3\ell_{3}}{27\ell_{3}}$,
 } &  & \tabularnewline
{\small$\plist[18\ell_{3}]{3\ell_{3}}{6\ell_{3}}{5\ell_{3}}{25\ell_{3}}$,
 } &  & \tabularnewline
{\small$\plist[21\ell_{3}]{3\ell_{3}}{9\ell_{3}}{8\ell_{3}}{28\ell_{3}}$.
 } &  & \tabularnewline
\end{tabular}}{\small}\\
{\small\par}
\end{enumerate}
\end{thm}

Notice that inside every expression in the above list at least one
of the exponents $i$, $j$, $s$, and $t$ is odd when $\ell,\ell_{1},\ell_{2},$
or $\ell_{3}$, respectively, are odd. We have\LetLtxMacro\itemold\item
\renewcommand{\item}{\itemindent0cm\itemold}
\begin{thm}
\label{th:OLodd} If\/ $k$ is odd, only the trivial overlaps occurs.
\end{thm}

\begin{proof}
Assume on the contrary that $(i,j\mid s,t)\in\cO(\bC,k)$ is reduced
with $k$ odd. Then $o:=(2i,2j\mid2s,2t)\in\cO(\bC,2k)$, where we
refer to a $2k$-th root of unity $\psi$ with $\psi^{2}=\phi$. Normalizing
$o$ would require at most one transformation of the form $\kappa_{i}:x\mapsto2k-x$,
which results in an even entry. The same holds for permutations when
applying to $o$. Thus, the normalized form of $o$ has four even
entries.

By Theorem~\ref{th.main} we have $2k=6\ell$. Since $k$ is odd,
$\ell$ is also odd. As we have noted, in this case, the list in Theorem~\ref{th.main}
shows no instance having four even entries. Thus there cannot be such
$(i,j\mid s,t)$ in $\cO(\bC,k)$.
\end{proof}

\section{The Proof}

\label{sec:The-Proof-Main-Theorem}

The working of the case (1) in the Theorem~\ref{th.main} has already
been verified in Example~\ref{ex:overlap}. It is also in \cite{KeK96}.
All others can be verified by similar methods based on the corresponding
cyclotomic polynomials. The main part of this section is to prove
that there are no more. We give some

\subsection{Preparations}

It will turn out to be convenient to make the following substitution
\[
a:=\h[s+j],~b:=\h[s-j],~c:=\h[t+i],~d:=\h[t-i],
\]
therefore
\begin{equation}
i=c-d,~j=a-b,~s=a+b,~t=c+d.\label{eq:uvst}
\end{equation}

We collect some easy consequences.
\begin{lem}
\label{l:le-abcd} Our assumptions on $i,j,s,t$ give
\begin{enumerate}
\item ${\displaystyle 0\le b<a<c\,,~0\le b<d<c\textnormal{\,,~\mbox{and}~}0\le b<d<\frac{k}{4}}$;
\item \textbf{${\displaystyle \frac{2\pi}{k}b<\frac{2\pi}{k}d\leq\pi-\frac{2\pi}{k}c<\pi-\frac{2\pi}{k}a}$;}
\item \textbf{${\displaystyle \pi-\frac{2\pi}{k}a+\frac{2\pi}{k}b<\pi-\frac{2\pi}{k}c+\frac{2\pi}{k}d}$.}
\end{enumerate}
\end{lem}

\begin{proof}
(1) follows easily from $0<i<j\le s<t\le\frac{k}{2}$ and $j+s<i+t$.
See Lemma~\ref{l:sumless}.

(2) Only $\frac{2\pi}{k}d\leq\pi-\frac{2\pi}{k}c$ needs explanation.
We have
\[
(\pi-\frac{2\pi}{k}c)-\frac{2\pi}{k}d=\pi-\frac{2\pi}{k}(c+d)=\pi-\frac{2\pi}{k}t\geq0.
\]

(3) ${\displaystyle \pi-\frac{2\pi}{k}a+\frac{2\pi}{k}b=\pi-\frac{2\pi}{k}j<\pi-\frac{2\pi}{k}i=\pi-\frac{2\pi}{k}c+\frac{2\pi}{k}d}$.
\end{proof}
By Lemma~\ref{l:polardec} the principle argument of $(\vpc[s]-1)^{-1}(\vpc[t]-1)$
is $(t-s)\pi/k$ while that of $(\vpc[i]-1)^{-1}(\vpc[j]-1)$ is $(j-i)\pi/k$.
Thus, $\vpc[\omega]=\exp\bigl(((t-s)-(j-i))\pi\bm{i}/k\bigr)$ and
so
\begin{equation}
\o=\frac{(t+i)-(s+j)}{2}=c-a.\label{omega}
\end{equation}

\begin{rem}
It turns out that $\omega$ is an integer. If $p\not\in\cQ_{k}$ and
$F=\text{GF}(q)$, where $q$ is a power of $p$ with $k\mid(q-1)$,
we have $\cO(F,k)=\cO(\mathbb{C},k)$. Moreover, if $\varphi\in F$
is a primitive $k$-th root of unity, and $(i,j\mid s,t)\in\cO(F,k)$,
then $\omega=\frac{(t+i)-(s+j)}{2}$ satisfies $\varphi^{\omega}\frac{1-\varphi^{j}}{1-\varphi^{i}}=\frac{1-\varphi^{t}}{1-\varphi^{s}}$.
Thus we know exactly how to compute $\omega$ from $i,j,s,t$. This
is not the case when $p\in\cQ_{k}$. Examples~\ref{ex:p13k7} and~\ref{ex:p11k12}
show such situations.
\end{rem}

Now, we expand (\ref{E.olp}) to obtain (see also \eqref{e:poly1})
\begin{equation}
1+\vpc[t+i]-\vpc[t]-\vpc[i]=\vpc[\omega]+\vpc[\omega+j+s]-\vpc[\omega+s]-\vpc[\omega+j].\label{E.expand_ori}
\end{equation}
Using $a,b,c,d$ and rearranging, we have
\begin{equation}
1+\vpc[2c]+\vpc[c+b]+\vpc[c-b]=\vpc[c+d]+\vpc[c-d]+\vpc[c-a]+\vpc[c+a].\label{E.expan}
\end{equation}
Multiply $\vpc[-c]$ to (\ref{E.expan}) and rearrange again to get
\[
\vpc[-a]+\vpc[a]+\vpc[-d]+\vpc[d]=\vpc[-b]+\vpc[b]+\vpc[-c]+\vpc[c].
\]
As $\vpc[-x]$ is the complex conjugate of $\vpc[x]$ for all $x$,
after dividing the last equation by $2$, we obtain
\begin{equation}
\cos\frac{2\pi}{k}a+\cos\frac{2\pi}{k}d=\cos\frac{2\pi}{k}b+\cos\frac{2\pi}{k}c.\label{eq:cosabc}
\end{equation}

This suggests that we shall be able to apply the following theorem
of Conway and Jones.
\begin{thm}[{{{{{{\cite[Theorem~7]{ConwayJ76}}}}}}}]
\label{Th.CJ7} Suppose we have at most four distinct rational multiples
of $\pi$ lying strictly between $0$ and $\frac{\pi}{2}$ for which
some rational linear combination of their cosines is rational but
no proper subset has this property. Then the appropriate linear combination
is proportional to one from the following list:
\begin{align}
\cos\frac{\pi}{3} & =\h[1],\label{cos:c0}\\
-\cos\theta+\cos\left(\frac{\pi}{3}-\theta\right)+\cos\left(\frac{\pi}{3}+\theta\right) & =0,\quad\left(0<\theta<\frac{\pi}{6}\right),\label{cos:c1}\\
\cos\frac{\pi}{5}-\cos\frac{2\pi}{5} & =\h[1],\label{cos:c2}\\
\cos\frac{\pi}{5}-\cos\frac{\pi}{15}+\cos\frac{4\pi}{15} & =\h[1],\label{cos:c4}\\
-\cos\frac{2\pi}{5}+\cos\frac{2\pi}{15}-\cos\frac{7\pi}{15} & =\h[1],\label{cos:c44}\\
\cos\frac{\pi}{7}-\cos\frac{2\pi}{7}+\cos\frac{3\pi}{7} & =\h[1].\label{cos:c5}
\end{align}
All other sums with four cosines equal $\h[1]$.
\end{thm}

Since (\ref{eq:cosabc}) has exactly four terms, we listed only the
relevant identities in the above theorem. In order to match \eqref{eq:cosabc}
with the equations in this theorem, we will have to rearrange the
equations so that all arguments to the cosine are in the range $[0,\pi/2]$,
and all terms nonnegative. For \eqref{eq:cosabc}, there are three
cases to consider.
\begin{lem}
\label{l:cosabc}
\begin{enumerate}
\item If $c\le\frac{k}{4}$, we just keep terms in \eqref{eq:cosabc} as
they are.
\item If $a\ge\frac{k}{4}$ and $c\ge\frac{k}{4}$, then \eqref{eq:cosabc}
must be transformed into
\begin{equation}
\cos\left(\pi-\frac{2\pi}{k}c\right)+\cos\frac{2\pi}{k}d=\cos\frac{2\pi}{k}b+\cos\left(\pi-\frac{2\pi}{k}a\right).\label{eq:cosabc-var}
\end{equation}
\item If $a<\frac{k}{4}\le c$, then \eqref{eq:cosabc} must be transformed
into
\begin{equation}
\cos\frac{2\pi}{k}a+\cos\frac{2\pi}{k}d+\cos\left(\pi-\frac{2\pi}{k}c\right)=\cos\frac{2\pi}{k}b.\label{eq:cosadc}
\end{equation}
\end{enumerate}
\end{lem}

\begin{proof}
By Lemma~\ref{l:le-abcd}\,(1), if $c\le\frac{k}{4}$, then all
terms in \eqref{eq:cosabc} are nonnegative; and, if $c>\frac{k}{4}$,
then according to $a<\frac{k}{4}$ and $a\ge\frac{k}{4}$, the other
two instances follow.
\end{proof}
With Lemma~\ref{l:cosabc}, we understand that we have to investigate
an expression of the form
\begin{equation}
e_{1}\cos\alpha_{1}+e_{2}\cos\alpha_{2}+e_{3}\cos\alpha_{3}+e_{4}\cos\alpha_{4}=0,\label{eq:cos-4}
\end{equation}
where $e_{i}\in\{1,-1\}$ and $\alpha_{i}$ are rational multiples
of $\pi$ in the range $[0,\pi/2]$. By Theorem~\ref{Th.CJ7}, equation
\eqref{eq:cos-4} must have subsums matching with equations \eqref{cos:c0}--\eqref{cos:c5}.
Notice that besides these possibilities there are always the trivial
equations
\[
\cos0=1\text{ }\mbox{ and }\text{ }\cos\h[\pi]=0,
\]
which can be used to fill up to four terms. By Lemma~\ref{l:le-abcd}(1),
$b$ is the single smallest value, thus $\cos0$ can occure at most
once in a sum.

We now go through the cases that can occur according to Theorem~\ref{Th.CJ7}.

\subsection{Cases (\ref{cos:c0}) and (\ref{cos:c1})}

We combine the first two equations of Theorem~\ref{Th.CJ7} to obtain
\begin{equation}
\cos\left(\frac{\pi}{3}-\theta\right)+\cos\left(\frac{\pi}{3}+\theta\right)=\cos\theta+\cos\frac{\pi}{2},~\mbox{where \ensuremath{0\le\theta<\frac{\pi}{6}}},\label{e:two.two}
\end{equation}
or
\begin{equation}
\cos\frac{\pi}{2}+\cos\left(\frac{\pi}{3}-\theta\right)+\cos\left(\frac{\pi}{3}+\theta\right)=\cos\theta,~\mbox{where \ensuremath{0\le\theta<\frac{\pi}{6}}}.\label{e:three.one}
\end{equation}
Notice that (\ref{cos:c0}) corresponds to the case $\theta=0$. Furthermore,
we have put in trivial terms to fill up to four. We also have
\[
\theta\le\frac{\pi}{3}-\theta\le\frac{\pi}{3}+\theta<\frac{\pi}{2}.
\]

By Lemma~\ref{l:le-abcd}(b), $\frac{2\pi}{k}b$ is the smallest
value in all three cases of Lemma~\ref{l:cosabc}. Thus, $\theta=\frac{2\pi}{k}b$,
or $b=\frac{k\theta}{2\pi}$. For convenience, we define $\ell=\frac{k}{6}$.
This gives us $0\leq b<\frac{k}{12}=\frac{\ell}{2}$. It will turn
out that $\ell$ is actually an integer.

Suppose that $c\le\frac{k}{4}$. By Lemmas~\ref{l:le-abcd}\,(1),
$c$ is the largest value. Therefore, \ref{l:cosabc} implies $\frac{2\pi}{k}c=\frac{\pi}{2}$,
and so $c=\frac{k}{4}$. Hence we can assume that $c\geq\frac{k}{4}$.
Now there are the following three possibilities:
\[
b<a\leq d<\frac{k}{4}\leq c,~b<d<a<\frac{k}{4}\leq c\text{ }\text{ or }\text{ }b<d<\frac{k}{4}\leq a<c.
\]

\noindent\textbf{(1)} $b<a\leq d<\frac{k}{4}\leq c$.

In this case, (\ref{eq:cosadc}) applies, and we have
\[
\cos\frac{2\pi}{k}a+\cos\frac{2\pi}{k}d+\cos\frac{2\pi}{k}\left(\frac{k}{2}-c\right)=\cos\frac{2\pi}{k}b.
\]
Matching up the arguments of this with those in (\ref{e:three.one}),
we have
\begin{align*}
\frac{2\pi}{k}\left(\frac{k}{2}-c\right)=\frac{\pi}{2},\quad\frac{2\pi}{k}d=\frac{\pi}{3}+\theta\text{ }\text{ and }\text{ }\frac{2\pi}{k}a=\frac{\pi}{3}-\theta.
\end{align*}
Therefore, remembering that $\theta=\frac{2\pi}{k}b$, one gets
\begin{align*}
c=\frac{k}{4}=\frac{3}{2}\ell,~d=\frac{k}{6}+b=\ell+b\text{ }\text{ and }\text{ }a=\frac{k}{6}-b=\ell-b.
\end{align*}
By \eqref{eq:uvst}, we get
\[
i=c-d=\frac{1}{2}\ell-b,\ j=a-b=\ell-2b
\]
and
\[
\ s=a+b=\ell,\ t=c+d=\frac{5}{2}\ell+b.
\]
Put $b'=\frac{1}{2}\ell-b$ to get
\[
(i,j\mid s,t)=(b',2b'\mid\ell,3\ell-b')\in\cO,~0<b'\leq\frac{1}{2}\ell.
\]
Now $\o=c-a=\frac{3\ell-b'+b'-(\ell+2b')}{2}=\ell-b'$, and so
\begin{equation}
\pp{\ell-b'}(2b',b',3\ell-b',\ell).\label{c3-1}
\end{equation}

\noindent\textbf{(2)} $b<d<a<\frac{k}{4}\leq c$.

Again, (\ref{eq:cosadc}) is used to match up with (\ref{e:three.one}),
and
\[
\frac{2\pi}{k}\left(\frac{k}{2}-c\right)=\frac{\pi}{2},~\frac{2\pi}{k}a=\frac{\pi}{3}+\theta\text{ }\text{ and }\text{ }\frac{2\pi}{k}d=\frac{\pi}{3}-\theta.
\]
Therefore,
\[
c=\frac{k}{4}=\frac{3}{2}\ell,~a=\frac{k}{6}+b=\ell+b\text{ }\text{ and }\text{ }d=\frac{k}{6}-b=\ell-b.
\]
By \eqref{eq:uvst}, we get
\[
i=c-d=\frac{1}{2}\ell+b,~j=a-b=\ell,~s=a+b=\ell+2b\text{ }\text{ and }\text{ }t=c+d=\frac{5}{2}\ell-b.
\]
Put $b''=\frac{1}{2}\ell+b$ to get
\[
(i,j\mid s,t)=(b'',\ell\mid2b'',3\ell-b'')\in\cO,~\frac{1}{2}\ell\leq b''<\ell.
\]
Notice that $b'=\frac{\ell}{2}$ in the previous case and $b''=\frac{\ell}{2}$
here give the same $(i,j\mid s,t)$.

Now $\o=c-a=\frac{3\ell-b''+b''-(\ell+2b'')}{2}=\ell-b''$ and
\begin{equation}
\pp{\ell-b''}(\ell,b'',3\ell-b'',2b'').\label{c3-2}
\end{equation}

\noindent\textbf{(3)} $b<d<\frac{k}{4}\leq a<c$.

In this case, (\ref{eq:cosabc-var}) applies, and we have
\[
\cos\frac{2\pi}{k}d+\cos\frac{2\pi}{k}\left(\frac{k}{2}-c\right)=\cos\frac{2\pi}{k}b+\cos\frac{2\pi}{k}\left(\frac{k}{2}-a\right).
\]
Matching up the arguments of this with those in (\ref{e:two.two}),
we have
\begin{align*}
\frac{2\pi}{k}d=\frac{\pi}{3}-\theta,~\frac{2\pi}{k}\left(\frac{k}{2}-c\right)=\frac{\pi}{3}+\theta\text{ }\text{ and }\text{ }\frac{2\pi}{k}\left(\frac{k}{2}-a\right)=\frac{\pi}{2}.
\end{align*}
Therefore,
\begin{align*}
d=\frac{k}{6}-\frac{k\theta}{2\pi}=\ell-b,~c=\frac{k}{3}-\frac{k\theta}{2\pi}=2\ell-b\text{ }\text{ and }\text{ }a=\frac{k}{4}=\frac{3}{2}\ell.
\end{align*}
By \eqref{eq:uvst}, we get
\[
i=c-d=\ell,\ j=a-b=\frac{3}{2}\ell-b,
\]
and
\[
s=a+b=\frac{3}{2}\ell+b,\ t=c+d=3\ell-2b.
\]
Put $b'''=\frac{3}{2}\ell-b$ to obtain
\[
(i,j\mid s,t)=(\ell,b'''\mid3\ell-b''',2b''')\in\cO,~\ell<b'''\leq\frac{3}{2}\ell.
\]
As $\o=c-a=\frac{2b'''+\ell-(b'''+3\ell-b''')}{2}=b'''-\ell$, we
find
\[
\pp{b'''-\ell}(b''',\ell,2b''',3\ell-b'''),
\]
or equivalently,
\begin{equation}
\pp{\ell-b'''}(\ell,b''',3\ell-b''',2b''').\label{c3-3}
\end{equation}

Note that $b',b'',b'''$ must be integers, as they are equal to one
of $i$ or $j$. Likewise, $\ell$ is an integer. Note also that the
ranges fit. Putting (\ref{c3-1}), (\ref{c3-2}), and (\ref{c3-3})
together, we therefore obtain
\[
\pp{\ell-u}(\ell,u,3\ell-u,2u),~1\leq u\leq\left\lfloor \frac{k}{4}\right\rfloor ,\;u\not=\frac{k}{6}.
\]
This is $\cO_{1}$.
\begin{rem}
In (1) and (2), when $\ell$ is even,
\[
\left(\frac{1}{2}\ell+b,\ell\Bigm|\ell+2b,\frac{5}{2}\ell-b\right)=\left(\frac{1}{2}\ell-b,\ell-2b\Bigm|\ell,\frac{5}{2}\ell+b\right)
\]
if and only if $b=0$. Obviously, the expression in~(3) cannot be
equal to one from (1) or~(2).
\end{rem}

\subsection{Case (\ref{cos:c2})}

We can use $\cos0=1$ or $\cos\frac{\pi}{2}=0$, to fill up (\ref{cos:c2})
to four terms. By Lemma~\ref{l:le-abcd}\,(1), $b$ is the smallest
values that appear among the arguments of $\cos$, if $\cos0$ is
used, we have $b=0$; otherwise, $b\not=0$.

In the case when $b\not=0$, (\ref{cos:c2}) from Theorem~\ref{Th.CJ7}
reads
\begin{equation}
\cos\frac{\pi}{5}+\cos\frac{\pi}{2}=\cos\frac{2\pi}{5}+\cos\frac{\pi}{3}.\label{eq:2.5-2:7}
\end{equation}
To simplify notation we let $k=60\ell_{3}$.

If $c\le\frac{k}{4}=15\ell_{3}$, we can match (\ref{eq:2.5-2:7})
with (\ref{eq:cosabc}). Using Lemma~\ref{l:le-abcd}(1), we obtain
\[
b=\frac{k}{10}=6\ell_{3},~c=\frac{k}{4}=15\ell_{3},\text{ }\text{ and }\text{ }\{a,d\}=\{10\ell_{3},12\ell_{3}\}.
\]
From this we get
\[
(i,j\mid s,t)=(3\ell_{3},4\ell_{3}\mid16\ell_{3},27\ell_{3})\mbox{ or }(5\ell_{3},6\ell_{3}\mid18\ell_{3},25\ell_{3}).
\]

If $c>\frac{k}{4}=15\ell_{3}$ we match (\ref{eq:2.5-2:7}) with (\ref{eq:cosabc-var}).
Using Lemma~\ref{l:le-abcd}\,(2), we obtain
\[
b=6\ell_{3},~d=10\ell_{3},~\frac{k}{2}-a=\frac{k}{4}\Longrightarrow a=15\ell_{3}\text{ }\text{ and }\text{ }c=18\ell_{3},
\]
leading to $(i,j\mid s,t)=(8\ell_{3},9\ell_{3}\mid21\ell_{3},28\ell_{3})$.

We have all possibilities in $\cO_{60}$. Notice that $\ell_{3}$
must be an integer, as the entries inside our three quadruples are
relatively prime.

Now we do the case $b=0$. In this case $\cos0$ has to be added to
(\ref{cos:c2}) from Theorem~\ref{Th.CJ7} to make it an equation
with four terms, which then reads
\begin{equation}
\cos\frac{\pi}{5}+\cos\frac{\pi}{3}=\cos\frac{2\pi}{5}+\cos0.\label{eq:new28}
\end{equation}
To simplify notation we let $k=30\ell_{1}$.

If $c\leq\frac{k}{4}$, we can match (\ref{eq:new28}) with (\ref{eq:cosabc}).
Using Lemma~\ref{l:le-abcd}(1), we obtain
\[
b=0\ell_{1},~c=6\ell_{1}\text{ }\text{ and }\text{ }\{a,d\}=\{3\ell_{1},5\ell_{1}\}.
\]
From this we get $(i,j\mid s,t)=(\ell_{1},3\ell_{1}\mid3\ell_{1},11\ell_{1})\mbox{ or }(3\ell_{1},5\ell_{1}\mid5\ell_{1},9\ell_{1})$.

If $c>\frac{k}{4}$ we can match (\ref{eq:new28}) with (\ref{eq:cosabc-var}).
Using Lemma~\ref{l:le-abcd}(2) we obtain
\[
b=0\ell_{1},~d=3\ell_{1},~a=9\ell_{1}\text{ }\text{ and }\text{ }c=10\ell_{1},
\]
leading to $(i,j\mid s,t)=(7\ell_{1},9\ell_{1}\mid9\ell_{1},13\ell_{1})$.

We have obtained the first three possibilities in $\cO_{30}$.

\subsection{Case (\ref{cos:c4})}

Slightly rewriting the equation (\ref{cos:c4}) from Theorem~\ref{Th.CJ7}
reads
\begin{equation}
\cos\frac{\pi}{5}+\cos\frac{4\pi}{15}=\cos\frac{\pi}{15}+\cos\frac{\pi}{3}.\label{eq:2.5-2:8}
\end{equation}
To simplify notation again we let $k=30\ell_{1}$.

If $c\leq\frac{k}{4}$, we can match (\ref{eq:2.5-2:8}) with (\ref{eq:cosabc}).
Using Lemma~\ref{l:le-abcd}(1) we obtain
\[
b=\ell_{1},~c=5\ell_{1}\text{ }\text{ and }\text{ }\{a,d\}=\{3\ell_{1},4\ell_{1}\}.
\]
From this we get $(i,j\mid s,t)=(\ell_{1},2\ell_{1}\mid4\ell_{1},9\ell_{1})\mbox{ or }(2\ell_{1},3\ell_{1}\mid5\ell_{1},8\ell_{1})$.

If $c>\frac{k}{4}$ we can match (\ref{eq:2.5-2:8}) with (\ref{eq:cosabc-var}).
Using Lemma~\ref{l:le-abcd}(2) we obtain
\[
b=\ell_{1},~d=3\ell_{1},~a=10\ell_{1}\text{ }\text{ and }\text{ }c=11\ell_{1},
\]
leading to $(i,j\mid s,t)=(8\ell_{1},9\ell_{1}\mid11\ell_{1},14\ell_{1})$.

These make the second block of three in $\cO_{30}$.

\subsection{Case (\ref{cos:c44})}

Slightly rewriting the equation (\ref{cos:c44}) from Theorem~\ref{Th.CJ7}
reads
\begin{equation}
\cos\frac{2\pi}{15}=\cos\frac{5\pi}{15}+\cos\frac{6\pi}{15}+\cos\frac{7\pi}{15}.\label{eq:2.5-2:9}
\end{equation}
For simple notation we stay with $k=30\ell_{1}$.

Now, the equation (\ref{eq:2.5-2:9}) can only match with (\ref{eq:cosadc}),
thus $c\geq\frac{k}{4}$ and $a<\frac{k}{4}$.

First matching $\frac{2\pi b}{k}$ with $\frac{2\pi}{15}$, we get
$b=2\ell_{1}$. From Lemma~\ref{l:le-abcd}\,(2), we have $d\leq\frac{k}{2}-c$.
Thus, there are three possibilities.\par
\hspace*{-1.5cm}\begin{minipage}[b]{1.085\hsize}
\begin{itemize}
\item If $d\leq\frac{k}{2}-c<a$, then $a=7\ell_{1}$, $15\ell_{1}-c=6\ell_{1}$
and $d=5\ell_{1}$. In this case, the element in $\cO_{30}$ (the
third block) is $(4\ell_{1},5\ell_{1}\mid9\ell_{1},14\ell_{1})$.
\item If $d<a<\frac{k}{2}-c$, then $15\ell_{1}-c=7\ell_{1}$, $a=6\ell_{1}$
and $d=5\ell_{1}$. In this case, the element in $\cO_{30}$ (the
third block) is $(3\ell_{1},4\ell_{1}\mid8\ell_{1},13\ell_{1})$.
\item If $a<d\leq\frac{k}{2}-c$, then $15\ell_{1}-c=7\ell_{1}$, $d=6\ell_{1}$
and $a=5\ell_{1}$. In this case, the element in $\cO_{30}$ (the
third block) is $(2\ell_{1},3\ell_{1}\mid7\ell_{1},14\ell_{1})$.
\end{itemize}
\end{minipage}\\
As before we emphasize the point that the entries inside all nine
quadruples involving $\ell_{1}$ are relatively prime. Thus~$\ell_{1}$
turns out to be an integer.

\subsection{Case (\ref{cos:c5})}

Slightly rewriting the equation (\ref{cos:c5}) in Theorem~\ref{Th.CJ7}
we get
\begin{equation}
\cos\frac{\pi}{7}+\cos\frac{3\pi}{7}=\cos\frac{2\pi}{7}+\cos\frac{\pi}{3}.\label{eq:2.5-c5}
\end{equation}
To simplify notation here, we let $k=42\ell_{2}$.

If $c\le\frac{k}{4}$, we can match the terms of (\ref{eq:2.5-c5})
with \eqref{eq:cosabc} and obtain
\[
b=3\ell_{2},~c=9\ell_{2},~\{a,d\}=\{6\ell_{2},7\ell_{2}\}.
\]
From this we get
\[
(i,j\mid s,t)=(2\ell_{2},3\ell_{2}\mid9\ell_{2},16\ell_{2})\mbox{ or }(3\ell_{2},4\ell_{2}\mid10\ell_{2},15\ell_{2}).
\]

If $c>\frac{k}{4}$, we can use Lemma~\ref{l:le-abcd}(2) again to
match (\ref{eq:2.5-c5}) with (\ref{eq:cosabc-var}) and get
\[
b=3\ell_{2},~d=\frac{1}{7}k=6\ell_{2},~c=14\ell_{2}\text{ }\text{ and }\text{ }a=12\ell_{2}.
\]
Therefore we obtain one more solution $(i,j\mid s,t)=(8\ell_{2},9\ell_{2}\mid15\ell_{2},20\ell_{2})$
to form $\cO_{42}$.

Again, $\ell_{3}$ must be an integer.

We have now exhibited all possible nontrivial overlaps. There are
no more than those listed in the theorem as we have claimed.

\section{Triple Overlaps}

\label{sec:Triple-Overlaps}

We say that a nontrivial \textit{triple overlap} occurs if for some
$s_{i},t_{i}\in\bk$, $i=1,2,3$, $(s_{1},t_{1}\mid s_{2},t_{2})$,
$(s_{2},t_{2}\mid s_{3},t_{3})$, and $(s_{1},t_{1}\mid s_{3},t_{3})$
are in $\cO$. In this case, we write $(s_{1},t_{1}\mid s_{2},t_{2}\mid s_{3},t_{3})$
and call it a \textit{nontrivial} triple overlap.

We collect all nontrivial triple overlaps in the set $\cT$, i.e.,
\begin{align*}
(s_{1},t_{1} & \mid s_{2},t_{2}\mid s_{3},t_{3})\in\cT\\
 & \Longleftrightarrow(s_{1},t_{1}\mid s_{2},t_{2}),(s_{2},t_{2}\mid s_{3},t_{3}),(s_{1},t_{1}\mid s_{3},t_{3})\in\cO.
\end{align*}
We are interested in nontrivial overlaps only. So if $(s_{1},t_{1}\mid s_{2},t_{2}\mid s_{3},t_{3})$
is in $\cT$, we simply use the phrase
\begin{center}
  ``there is a triple overlap
$(s_{1},t_{1}\mid s_{2},t_{2}\mid s_{3},t_{3})$,''
\end{center} or
\begin{center}
  ``$(s_{1},t_{1}\mid s_{2},t_{2}\mid s_{3},t_{3})$
is a triple overlap'',
\end{center}and the like to refer to a nontrivial triple
overlap.

Recall that when $k$ is odd, no nontrivial overlaps can occur by
Theorem~\ref{th:OLodd}. We are only dealing with even~$k$. Thus,
we may assume that $s_{i}\leq\frac{k}{2}$ and $t_{i}\leq\frac{k}{2}$
for all $i$.

For a triple overlap $(s_{1},t_{1}\mid s_{2},t_{2}\mid s_{3},t_{3})\in\cT$,
we denote
\[
o_{1}=(s_{1},t_{1}\mid s_{2},t_{2}),\quad o_{2}=(s_{2},t_{2}\mid s_{3},t_{3}),\quad o_{3}=(s_{1},t_{1}\mid s_{3},t_{3}),
\]
and call them the \textit{constituents} of the triple overlap.

The following are easy consequences from the definition.
\begin{lem}
\label{l:TO-Trans} It holds that
\begin{align*}
(s_{1},t_{1}\mid s_{2},t_{2}\mid s_{3},t_{3})\in\cT & \Longleftrightarrow\hbox{}(s_{2},t_{2}\mid s_{1},t_{1}\mid s_{3},t_{3})\in\cT\\
 & \Longleftrightarrow\hbox{}(s_{3},t_{3}\mid s_{2},t_{2}\mid s_{1},t_{1})\in\cT.
\end{align*}
\end{lem}

\begin{lem}
\label{l:TO-switch} It holds that
\[
(s_{1},t_{1}\mid s_{2},t_{2}\mid s_{3},t_{3})\in\cT\Longleftrightarrow(t_{1},s_{1}\mid t_{2},s_{2}\mid t_{3},s_{3})\in\cT.
\]
\end{lem}

To obtain all the triple overlaps in $\cT$, we can therefore restrict
to the case when $s_{1}<s_{2}<s_{3}$ and $s_{1}<t_{1}$. More precisely,
if there is a triple overlap $(s_{1}',t_{1}'\mid s_{2}',t_{2}'\mid s_{3}',t_{3}')$,
then there is also a triple overlap $(s_{1},t_{1}\mid s_{2},t_{2}\mid s_{3},t_{3})$
with $s_{1}<s_{2}<s_{3}$ and $s_{1}<t_{1}$, and this one is referred
to as \textit{normalized}.
\begin{lem}
\label{l:TO-normal} Let $(s_{1},t_{1}\mid s_{2},t_{2}\mid s_{3},t_{3})\in\cT$
be a normalized triple overlap. Then $s_{i}<t_{i}$ and $t_{1}<t_{2}<t_{3}$.
\end{lem}

\begin{proof}
Lemma~\ref{l:order-uv} implies $s_{i}<t_{i}$ for all $i$ and $t_{1}<t_{2}<t_{3}$.
\end{proof}
\begin{thm}
\label{t:tripolp} Let $T=(s_{1},t_{1}\mid s_{2},t_{2}\mid s_{3},t_{3})\in\cT$
be a normalized triple overlap. Then $k$ is divisible by $30$, i.e.,
$k=30\ell$, and $T$ is one of the following
\begin{alignat*}{2}
T_{1}: & \ttolp[3]5596{12}, & \qquad T_{4}: & \ttolp[2]5384{13},\\
T_{2}: & \ttolp 2495{14}, & T_{5}: & \ttolp[4]58{11}9{14}.\\
T_{3}: & \ttolp[2]3587{14}, &  & \negthickspace
\end{alignat*}
Any other triple overlap in $\cT$ can be obtained from these by applying
the operations from Lemmas~\ref{l:TO-Trans} and~\ref{l:TO-switch}.
\end{thm}

By inspection, we derive from this an immediate consequence.
\begin{cor}
\label{cor:quadOL} There do not exist nontrivial ``quadruple''
overlaps.
\end{cor}

\section{Proof of Theorem~\ref{t:tripolp}}

By Lemmas~\ref{l:TO-Trans}, \ref{l:TO-switch} and~\ref{l:TO-normal}
every triple overlap $(s_{1}',t_{1}'\mid s_{2}',t_{2}'\mid s_{3}',t_{3}')$
in $\cT$ can be normalized into a triple overlap $(s_{1},t_{1}\mid s_{2},t_{2}\mid s_{3},t_{3})$
with $s_{1}<s_{2}<s_{3}$, $t_{1}<t_{2}<t_{3}$, and $s_{i}<t_{i}$
for $i=1,2,3$. This accounts for the last statement.

Throughout this proof, overlaps and triple overlaps are not necessarily
normalized. Thus, when doing inspection below, this fact has to be
taken into account. Let $(s_{1},t_{1}\mid s_{2},t_{2}\mid s_{3},t_{3})\in\cT$.
Whenever we are dealing with overlaps from $\cO_{30}$, or $\cO_{42}$,
or $\cO_{60}$, we assume that $k=30\ell_{1}$, $k=42\ell_{2}$, $k=60\ell_{3}$,
respectively.

The cases are organized by the number of constituents inside $\cO_{1}$.

\subsection{Triple overlaps with no constituents in $\protect\cO_{1}$}

\noindent\textbf{(1)} Assume that there is a constituent, $o_{1}$
say, in $\cO_{60}$. Then using Theorem~\ref{th.main} one easily
sees that the others are not in~$\cO_{60}$.

If another constituent is to be in $\cO_{30}$, then $\ell_{1}=2\ell_{3}$
and $o_{1}$ must have one instance of the form ${(2m,2m'\mid.\,,\,.\,)}$.
However, this is not the case according to Theorem~\ref{th.main}.
Similarly, if another constituent is to be in $\cO_{42}$, then $\ell_{2}$
is a multiple of $10$ and so two entries of $o_{1}$ must be multiples
of $10$, which is not the case either.

Therefore, no constituents can be in $\cO_{60}$.

\noindent\textbf{(2)} Assume that there is a constituent, $o_{1}$
say, in $\cO_{42}$. Then from the list in Theorem~\ref{th.main}
one easily sees that the others are not in~$\cO_{42}$. Hence, other
constituents $o_{2}$ and $o_{3}$ have to be in $\cO_{30}$. But
then $\ell_{1}$ is a multiple of $7$ and so two entries of $o_{2}$,
say, must be multiples of $7$, which does not happen.

Therefore, no constituents can be in $\cO_{42}$.

\noindent\textbf{(3)} We are left with the case that there are two
constituents in $\cO_{30}$. A tedious inspection (see the Remark~\ref{rem:strat}
below) reveals two normalized triple overlaps:
\begin{equation}
(\ell_{1},2\ell_{1}\mid4\ell_{1},9\ell_{1}\mid5\ell_{1},14\ell_{1})\mbox{\quad and\quad}(2\ell_{1},5\ell_{1}\mid3\ell_{1},8\ell_{1}\mid4\ell_{1},13\ell_{1}),\label{eq:noO1}
\end{equation}
as well as one not normalized: $(4\ell_{1},5\ell_{1}\mid9\ell_{1},14\ell_{1}\mid8\ell_{1},11\ell_{1})$.
Here $k=30\ell_{1}$.

In each of the three triple overlaps found, there is one constituent
(in fact $o_{3}$) from $\cO_{1}$. Thus, these triple overlaps do
not meet the condition of the present case, and will show up again
in the sequel.
\begin{rem}
\label{rem:strat} Here we describe an efficient strategy to find
normalized triple overlaps from the list of overlaps.
\begin{itemize}
\item Start with a normalized overlap $(s_{1},t_{1}\mid s_{2},t_{2})$.
\item Look for normalized overlaps $(s_{2},v\mid s,t)$ and check \\
 ${\quad}$if $t_{2}=v$, then $(s_{1},t_{1}\mid s_{2},t_{2}\mid s,t)\in\cT$,\quad{}and\\
 ${\quad}$if $t_{2}=s$, then $(s_{1},t_{1}\mid s_{2},t_{2}\mid v,t)\in\cT$.
\end{itemize}
Notice that the last of the triple overlaps given above cannot be
found this way. Yet, it is spotted during the inspection. A normalized
one for it will be found later.
\end{rem}

\subsection{Triple overlaps with exactly one constituent in $\protect\cO_{1}$}

Putting $k=6\ell$, we can assume that $o_{1}=(i,\ell\mid2i,3\ell-i)$.
This means that $o_{2}$ or $o_{3}$ must contain $\ell$. Furthermore,
we have $o_{2},o_{3}\in\cO_{30}\cup\cO_{42}\cup\cO_{60}$.

In the cases $k=60\ell_{3}$ and $k=42\ell_{2}$ this implies that
$10\ell_{3}$ or $7\ell_{2}$, respectively, must occur as an entry
of an overlap. Yet there is no such overlap in the list of Theorem~\ref{th.main}.

In the case $k=30\ell_{1}$, the entry $5\ell_{1}$ must occur in
an overlap inside $\cO_{30}$. There are three normalized candidates:
\begin{equation}
(3\ell_{1},5\ell_{1}\mid5\ell_{1},9\ell_{1}),\ (2\ell_{1},3\ell_{1}\mid5\ell_{1},8\ell_{1}),\ (4\ell_{1},5\ell_{1}\mid9\ell_{1},14\ell_{1}).\label{eq:newO1}
\end{equation}
If we assume that the first entry is smaller than the second, then
three pairs of this form $(i,5\ell_{1})$, $(5\ell_{1},i)$, $(5\ell_{1},15\ell_{1}-i)$
can be found in the overlaps in $\cO_{1}$. Note that $5\ell_{1}<15\ell_{1}-i$,
as $i\le15\ell_{1}/2$.

From the first quadruple of \eqref{eq:newO1}, we get $i=3\ell_{1}$,
or $i=9\ell_{1}$ (too large), or $15\ell_{1}-i=9\ell_{1}$, hence
$i=3\ell_{1}$ or $6\ell_{1}$. Both yield
\[
(3\ell_{1},5\ell_{1}\mid5\ell_{1},9\ell_{1}\mid6\ell_{1},12\ell_{1}).
\]
Since this triple overlap has two constituents inside $\cO_{1}$,
it does not fit the condition we are considering, and will show up
again in the next case.

From the second quadruple of \eqref{eq:newO1}, we get $i=2\ell_{1}$,
or $i=8\ell_{1}$ (too large), or $15\ell_{1}-i=8\ell_{1}$, hence
(already normalized)\smallskip{}
 \\
 ${~\qquad}$$i=2\ell_{1}$, giving $(2\ell_{1},5\ell_{1}\mid3\ell_{1},8\ell_{1}\mid4\ell_{1},13\ell_{1})$,
\\
 and\\
 ${~\qquad}$$i=7\ell_{1}$, giving $(2\ell_{1},3\ell_{1}\mid5\ell_{1},8\ell_{1}\mid7\ell_{1},14\ell_{1})$.
\smallskip{}
 \\
 Notice that the first one here is the second in \eqref{eq:noO1},
and both triple overlaps here have one constituent in $\cO_{1}$ and
two constituents in~$\cO_{30}$.

From the third quadruple of \eqref{eq:newO1}, we have $i=4\ell_{1}$,
or $i=14\ell_{1}$ (too large), or $15\ell_{1}-i=14\ell_{1}$, hence\smallskip{}
 \\
 ${~\qquad}$$i=1\ell_{1}$, giving $(\ell_{1},2\ell_{1}\mid4\ell_{1},9\ell_{1}\mid5\ell_{1},14\ell_{1})$,
\\
 and\\
 ${~\qquad}$$i=4\ell_{1}$, giving $(4\ell_{1},5\ell_{1}\mid8\ell_{1},11\ell_{1}\mid9\ell_{1},14\ell_{1})$.
\smallskip{}
 \\
 The first one here is the first one in \eqref{eq:noO1} while the
second one is the nonnormalized triple overlap we had after \eqref{eq:noO1}.
Both of them have one constituent in $\cO_{1}$ and two constituents
in~$\cO_{30}$.

Note that we have now found all the triple overlaps listed in the
theorem, including $(3\ell_{1},5\ell_{1}\mid5\ell_{1},9\ell_{1}\mid6\ell_{1},12\ell_{1})$
which did not fit the current case condition. The following discussions
will only show this one, but not any new ones.

\subsection{Triple overlaps with at least two constituents in $\protect\cO_{1}$}

Assume $k=6\ell$ and start out with a triple overlap $o\in\cT$ with
at least two constituents from $\cO_{1}$ (may not be normalized nor
reduced, and the numbering may not be the actual order):
\[
o_{1}=(i,\ell\mid2i,3\ell-i)\text{\quad and\quad}o_{2}=(j,\ell\mid2j,3\ell-j),
\]
where $1\leq i\leq\left\lfloor \frac{k}{4}\right\rfloor =\left\lfloor \frac{3\ell}{2}\right\rfloor <2\ell$,
$1\leq j\le\left\lfloor \frac{3\ell}{2}\right\rfloor $, $i\not=\ell$,
$j\not=\ell$, and $i\not=j$. Thus, we have $2i\leq3\ell$, $2j\leq3\ell$,
and $i+j\leq3\ell$. There is no loss of generality in assuming that
$i<j$. Then from $i<j\leq\frac{3\ell}{2}$, we infer that $i+j<3\ell$.
Summarizing, we obtain the following four inequalities:
\[
i<j,\quad i<3\ell-i,\quad\ell<3\ell-i,\text{ }\text{ and }\text{ }i<3\ell-j.
\]
From these, we see that either $i$ or $\ell$ is the smallest among
all entries involved in the triple overlap.

In both cases, $3\ell-i$ is the largest entry in $o_{1}$, and must
appear in $o_{2}$. That is, $3\ell-i$ must be one of $j$, $\ell$,
$2j$ and $3\ell-j$. The first, second and fourth cases lead to contradictions
$i+j=3\ell$, $i=2\ell$ and $i=j$, respectively. The third case
makes $3\ell-i=2j$. Also, since $i<j$, we have $3\ell=i+2j<3j$,
and so $\ell<j$. Consequently, $i<\ell$.

Either $(i,\ell\mid2i,3\ell-i)=(i,\ell\mid2i,2j)$ or $(i,\ell\mid2i,3\ell-i)\sim(i,2i\mid\ell,2j)$
is the normalized form of~$o_{1}$. Assume the later is normalized.
In order to match it, we have to rearrange $o_{2}=(j,\ell\mid2j,3\ell-j)$
into the form $(j,2j\mid\ell,3\ell-j)$ or $(3\ell-j,2j\mid\ell,j)$.
Then either $\ell=j$ or $\ell=3\ell-j$, contradicting the fact that
$\ell<j\leq\lfloor\frac{3\ell}{2}\rfloor$.

Therefore, $o_{1}=(i,\ell\mid2i,2j)$ is in normalized form. And again,
one of the rearrangements $(j,2j\mid\ell,3\ell-j)$ or $(3\ell-j,2j\mid\ell,j)$
of $o_{2}$ matches $o_{1}$. If $2i=3\ell-j$, then, together with
$3\ell-i=2j$, we arrive at the contradiction $i=j$. Thus we are
left with $2i=j$. This yields $i=\frac{3}{5}\ell$ and $j=\frac{6}{5}\ell$.
Putting $k=30\ell_{1}$, we obtain the triple overlap $(3\ell_{1},5\ell_{1}\mid6\ell_{1},12\ell_{1}\mid5\ell_{1},9\ell_{1})$,
which already showed up earlier.

After rearranging into normalized form we obtain $o_{1}\in\cO_{30}$,
$o_{2}\in\cO_{1}$ and $o_{3}\in\cO_{1}$.

Finally, there are no triple overlaps with all three constituents
inside $\cO_{1}$. We are done with the proof.

\section{Acknowledgements}

The authors thank the reviewer for his/her careful reading and comments.

\begin{thebibliography}{1}
\bibitem{Clay88}J. R. Clay. \emph{Circular block designs from planar
near-rings.} Combinatorics '86 (Trento, 1986), 95--105, Ann. Discrete
Math., \textbf{37}, North-Holland, Amsterdam, 1988.

\bibitem{Clay92}J. R. Clay. Nearrings: Geneses and Applications.
Oxford Univ. Press, Oxford, 1992.

\bibitem{ConwayJ76}J. H. Conway and A. J. Jones. \emph{Trigonometric
diophantine equations}, Acta Arith. \textbf{30} (1976), 229--240.

\bibitem{IreRos} K. F. Ireland and M. I. Rosen. \emph{A classical
introduction to modern number theory}, $2^{\text{nd}}$ ed., Springer-Verlag,
Berlin-Heidelberg-New York, 1990.

\bibitem{KeK96}W.-F. Ke and H. Kiechle, \emph{Combinatorial properties
of ring generated circular planar nearrings}, J. Combin. Theory Ser.
A \textbf{73} (1996), 286--301.

\bibitem{KeK95}W.-F. Ke and H. Kiechle, \textit{On the solutions
of the equation $x^{m}+y^{m}-z^{m}=1$ in a finite field,} Proc. Amer.
Math. Soc. 128 (1995), 1331--1339.

\bibitem{KeK23-1}W.-F. Ke and H. Kiechle, \textit{Circularity in
Finite Fields and Solutions of the Equations $x^{m}+y^{m}-z^{m}=1$.}
Submitted; \href{https://arxiv.org/abs/2307.05586}{arXiv:2307.05586}.

\bibitem{Kiechle94}H. Kiechle, \textit{Points on Fermat curves over
finite fields,} in ``Proc. Conference on Finite Fields: Theory, Applications,
and Algorithms, Las Vegas, NV, 1993,'' Contemp. Math. 168 (1994),
181--183.

\bibitem{Modisett89} M. C. Modisett. \emph{A characterization of
the circularity of balanced incomplete block designs.} Utilitas Math.
\textbf{35} (1989), 83--94.

\end{thebibliography}
\end{document}